\documentclass{svmult}

\renewcommand{\email}[1]{\emailname: #1} 

\usepackage{mathptmx}       
\usepackage{helvet}         
\usepackage{courier}        

\usepackage{makeidx}         
\usepackage{graphicx}        
\usepackage[bottom]{footmisc}

\usepackage{latexsym}
\usepackage{amsmath}
\usepackage{amsfonts}
\usepackage{amssymb}
\usepackage{bm}

\usepackage{url}
\usepackage{algorithm}
\usepackage{algorithmic}
\usepackage[misc,geometry]{ifsym}

\spdefaulttheorem{assumption}{Assumption}{\upshape \bfseries}{\itshape}
\spdefaulttheorem{algo}{Algorithm}{\upshape \bfseries}{\itshape}

\newcommand{\rd}{{\mathrm{d}}}

\newcommand{\bbE}{{\mathbb{E}}}

\DeclareSymbolFont{bbold}{U}{bbold}{m}{n}
\DeclareSymbolFontAlphabet{\mathbbold}{bbold}

\providecommand{\argmax}{\operatorname*{argmax}}

\begin{document}

\title*{Modern Monte Carlo Variants for Uncertainty Quantification in Neutron Transport}

\author{Ivan G. Graham \and Matthew J. Parkinson \and Robert Scheichl}

\institute{
 Ivan G. \ Graham (\Letter) \and Matthew J. \ Parkinson \and Robert Scheichl
 \at University of Bath, Claverton Down, Bath, BA2 7AY, UK\\
 \email{i.g.graham@bath.ac.uk; m.parkinson@bath.ac.uk; r.scheichl@bath.ac.uk}
}

\maketitle

\index{Graham, I.G.}
\index{Parkinson, M.J.}
\index{Scheichl, R.}

\paragraph{Dedicated to Ian H.~Sloan on the occasion of his 80th birthday.}

\abstract{
We describe modern variants of Monte Carlo methods for Uncertainty Quantification (UQ) of the Neutron Transport Equation, when it is approximated by the discrete ordinates method with diamond differencing. 
We focus on the mono-energetic 1D slab geometry problem, with isotropic
scattering,  where the cross-sections are  log-normal correlated random fields of possibly  
low regularity.  The paper includes an outline of  novel theoretical results on the convergence of the discrete scheme, in the cases  of both spatially variable and 
random cross-sections. 
We also describe the theory and practice of algorithms for  quantifying  
the uncertainty of a linear functional of the
scalar flux, using  Monte Carlo and quasi-Monte Carlo
methods, and their multilevel variants. 
A hybrid iterative/direct
solver for computing each realisation of the functional is also presented. Numerical
experiments  show the effectiveness of the hybrid solver
and the gains that are possible through quasi-Monte Carlo sampling and
multilevel variance reduction. For the multilevel quasi-Monte Carlo
method, we observe gains in the computational $\varepsilon$-cost of 
up to 2 orders of magnitude over the standard Monte Carlo method,  and we explain this theoretically.  
Experiments on problems with 
up to  several thousand stochastic dimensions are included.
}


\section{Introduction}

In this paper we will consider the Neutron Transport equation (NTE), sometimes referred to as the Boltzmann transport equation. This is an integro-differential equation which models the flux of neutrons in a reactor. It has particular applications for nuclear reactor design, radiation shielding and astrophysics \cite{SaMc:82}. There are many potential sources of uncertainty in a nuclear reactor, such as the geometry, material composition and reactor wear. Here, we will consider the problem of random spatial variation in the coefficients (the \textit{cross-sections}) in the NTE, represented by correlated random fields with potentially low smoothness. Our aim is to understand how \textit{uncertainty in the cross-sections propagates through to (functionals of) the neutron flux}. This is the forward problem of Uncertainty Quantification.

We will quantify the uncertainty using Monte Carlo (MC) type methods, that is, by simulating a finite number of pseudo-random instances of the NTE and by averaging the outcome of those simulations to obtain statistics of quantities of interest. Each statistic can be interpreted as an expected value of some (possibly nonlinear) functional of the neutron flux with respect to the random cross-sections. The input random fields typically need to be parametrised with a significant number of random parameters leading to a problem of high-dimensional integration. MC methods are known to be particularly well-suited to this type of problem due to their dimension independent convergence rates.


However, convergence of the MC algorithm is slow and determined by $\sqrt{ \mathbb{V}(\cdot) / N  \ }$, where $\mathbb{V}(\cdot)$ is the variance of the quantity of interest  and $N$ is the number of samples. For this reason, research is focussed on improving the convergence, whilst retaining dimensional independence. Advances in MC methods can broadly be split into two main categories: improved sampling and variance reduction.
Improved sampling methods attempt to find samples that perform better than the  pseudo-random choice. Effectively, they aim to improve the $\sqrt{ 1/N \ }$ term in the error estimate. A major advance in sampling methods has come through the development of quasi-Monte Carlo (QMC) methods. 
Variance reduction methods, on the other hand, attempt to reduce the $\mathbb{V}(\cdot)$ term in the error estimate and thus reduce the number of samples needed for a desired accuracy. Multilevel Monte Carlo (MLMC) methods (initiated in \cite{He:01, Gi:08} and further developed in, e.g., \cite{GiWa:09,BaScZo:11,Cl:11,ChSc:13,KuScSl:15,Kuo:15,TeJa:15,HaAli:15}) fall into this category. A comprehensive review of MLMC can be found in \cite{Gi:15}.

The rigorous theory of all of the improvements outlined above requires  regularity properties of the solution, the    
  verification of  which can be a substantial task. 
 There are  a significant number of published papers on the regularity of parametric elliptic PDEs, 
in physical and  parameter space, as they arise, e.g., in flow in random models of  porous media \cite{ChSc:13, KuScSl:12, DiKuGi:14a,DiKuGi:14,GrKu:15,KuScSl:15,Kuo:15}. However, for the NTE,  
this  regularity question  
is almost untouched. 
Our complementary paper \cite{GrPaSc:17} contains a full regularity and error analysis of the discrete scheme for the NTE with spatially variable and random coefficients. Here we restrict to a summary of those results. 

The field of UQ has grown very quickly in recent years and its application to neutron transport theory is currently of considerable interest. There are a number of groups that already work on this problem, e.g.  \cite{AyEa:15, Fi:11, Gi6:13} and references therein. 
Up to now,   research has focussed on using the polynomial chaos expansion (PCE), which comes in two forms; the intrusive and non-intrusive approaches.  Both approaches expand the random flux in a weighted sum of orthogonal polynomials. The intrusive approach considers the expansion directly in the differential equation, which in turn requires a new solver (`intruding' on the original solver). In contrast, the non-intrusive approach attempts to estimate the coefficients of the PCE directly, by projecting onto the PCE basis cf. \cite[eq.(40)]{AyEa:15}. This means the original solver can be used as a `black box' as in MC methods. 
Both of the approaches then use quadrature to estimate the coefficients in the PCE. The main disadvantage of standard PCE is that typically the number of terms grow exponentially in the number of stochastic dimensions and in the order of the PCE, the so-called \textit{curse of dimensionality}.

Fichtl and Prinja \cite{Fi:11} were some of the first to numerically tackle the 1D slab geometry problem with random cross-sections. 
Gilli et al. \cite{Gi6:13} improved upon this work by using (adaptive) sparse grid ideas in the collocation method, to tackle the curse of dimensionality. Moreover, \cite{AyPaEa:14} constructed a hybrid PCE using a combination of Hermite and Legendre polynomials, observing superior convergence in comparison to the PCE with just Hermite polynomials. 
More recently \cite{AyEa:15} tackled the (time-independent) full criticality problem in three spatial, two angular and one energy variable. They consider a second expansion, the high-dimensional model representation (HDMR), which allows them to expand the response (e.g. functionals of the flux) in terms of low-dimensional subspaces of the stochastic variable. The PCE is used on the HDMR terms, each with their own basis and coefficients. We note however, that none of these papers provide any rigorous error or cost analysis.

The structure of this paper is as follows. In Section 2, we describe the model problem, a 1D slab geometry simplification of the Neutron Transport Equation with spatially varying and random cross-sections. We set out the discretisation of this equation and discuss two methods for solving the resultant linear systems; a direct and an iterative solver.  In Section 3, the basic elements of a fully-discrete error analysis of the discrete ordinates method with diamond differencing applied to the model problem are summarised. The full analysis will be given in \cite{GrPaSc:17}. In Section 4, we introduce a number of variations on the Monte Carlo method for quantifying uncertainty. This includes a summary of the theoretical computational costs for each method. Finally, Section 5 contains numerical results relating to the rest of the paper. We first present a hybrid solver that combines the benefits of both  
direct and iterative solvers. Its cost depends on the particular realisation of the cross-sections. Moreover, we present  simulations for the UQ problem for the different variants of the Monte Carlo methods, and compare the rates with those given by the  theory.

\section{The Model Problem}
\label{sec:ModelProblem}

The Neutron Transport Equation (NTE) is a physically derived balance equation, that models the angular flux $\psi(\vec{r},\Theta,E)$ of neutrons in a domain, where $\vec{r}$ is position, $\Theta$ is angle and $E$ is energy. Neutrons are modelled as non-interacting particles travelling along straight line paths with some energy $E$. They interact with the larger nuclei via absorption, scattering  and fission. The rates $\sigma_A$, $\sigma_S$ and $\sigma_F$ at which these events occur are called the \textit{absorption, scattering and fission cross-sections}, respectively. They can depend on the position $\vec{r}$ and the energy $E$ of the neutron. The scattering cross-sections also depend on the energy $E'$ after the scattering event, as well as on the angles $\Theta$ and $\Theta'$ before and after the event. 

The two main scenarios of interest in neutron transport are the so-called \textit{fixed source problem} and the \textit{criticality problem}. We will focus on the former, which concerns the transport of neutrons emanating from some fixed source term $f$. It has particular applications in radiation shielding. 
We will further simplify our model to the \textit{1D slab geometry case} by assuming
\begin{itemize}
\item no energy dependence;
\item dependence only on one spatial dimension and infinite extent of the domain in the other two dimensions;
\item no dependence of any cross-sections on angle;
\item no fission.
\end{itemize}
The resulting simplified model is an integro-differential equation for the angular flux $\psi(x,\mu)$ such that 
\begin{align} \label{eq:transport_det} 
\mu \frac{\rd \psi}{\rd x}(x,\mu) \ + \ \sigma(x) \psi(x,\mu)  &\ = \ \sigma_S(x) \phi(x) \ +  \ f(x) \ , \\
\label{eq:phi_det}
\text{where} \qquad \phi(x) &\ = \ \frac{1}{2} \int_{-1}^1 \psi(x,\mu') \ \rd \mu' \ ,
\end{align}
for any $x \in (0,1)$ and $\mu \in [-1,1]$, subject to the no in-flow boundary conditions
\begin{equation}\label{eq:bc} 
\psi(0,\mu) \ = \ 0, \ \mbox{ for } \ \mu > 0 \quad \text{and} \quad \psi(1, \mu) = 0, \ \mbox{ for } \ \mu < 0 \ .
\end{equation}
Here, the angular domain is reduced from $\mathbb{S}_2$ to the unit circle $\mathbb{S}_1$ and parametrised by the cosine $\mu \in [-1,1]$ of the angle. The equation degenerates at $\mu = 0$, i.e. for neutrons moving perpendicular to the $x$-direction. The coefficient function $\sigma(x)$ is the total cross-section given by $\sigma = \sigma_S + \sigma_A$.  
For more discussion on the NTE see \cite{DaLi:12, LeMi:84}.

\subsection{Uncertainty Quantification}

\label{subsec:UQ} 

An important problem in industry is to quantify the uncertainty in the fluxes due to uncertainties in the cross-sections. Most materials, in particular shielding materials such as concrete, are naturally heterogeneous or change their properties over time through wear. Moreover, the values of the cross-sections are taken from nuclear data libraries across the world and they can differ significantly between libraries \cite{LeLe:07}. This means there are large amounts of uncertainty on the coefficients, and this could have significant consequences on the system itself. 

To describe the random model, let $(\Omega,\mathcal{A}, \mathbb{P})$ be a probability space with $\omega \in \Omega$ denoting a random event from this space. Consider a (finite) set of partitions of the spatial domain, where on each subinterval we assume that $\sigma_S = \sigma_S(x,\omega)$ and $\sigma = \sigma(x,\omega)$ are two 
(possibly dependent or correlated) random fields. Then the angular flux 
and the scalar flux become random fields and the model problem \eqref{eq:transport_det}, \eqref{eq:phi_det} becomes 
\begin{align} \label{eq:transport} 
\mu \frac{\rd \psi}{\rd x}(x,\mu,\omega) \ + \ \sigma(x, \omega)  \psi(x,\mu,\omega)  \ = \ \sigma_S(x,\omega) \phi(x,\omega) \ + \ f(x) \ , \\
\label{eq:phi}
\text{where} \qquad \phi(x,\omega) \ = \ \int_{-1}^1 \psi(x,\mu',\omega) d \mu' 
\end{align}
and $\psi(\cdot,\cdot,\omega)$ satisfies the boundary conditions \eqref{eq:bc}. The set of equations 
 \eqref{eq:transport}, \eqref{eq:phi}, \eqref{eq:bc}  have to hold for almost all realisations $\omega \in \Omega$.

For simplicity, we restrict ourselves to deterministic $\sigma_A = \sigma_A(x)$ with 
\begin{equation}
\label{eq:sigA}
0 \ < \ \sigma_{A,\min} \ \leq \ \sigma_A(x) \ \leq \ \sigma_{A,\max}  \ < \  \infty \ , \quad \text{for all} \quad x\in [0,1] \ , 
\end{equation}
and assume a log-normal distribution for $\sigma_S(x,\omega)$. The total cross-section $\sigma(x,\omega)$ is then simply the log-normal random field with values $\sigma(x,\omega) = \sigma_S(x,\omega) + \sigma_A(x)$.
 In particular, we assume that $\log \sigma_S$ is a correlated zero mean Gaussian random field, with covariance function defined by
\begin{equation}
\label{eq:covfunc}
C_\nu (x,y) \ = \ \sigma_{var}^2 \frac{2^{1-\nu}}{\Gamma(\nu)} \left( 2 \sqrt{\nu} \frac{|x - y|}{\lambda_C}  \right)^\nu K_\nu \left( 2 \sqrt{\nu} \frac{|x - y|}{\lambda_C} \right) \ .
\end{equation}
This class of covariances is called the Mat{\'e}rn class. It is parametrised by the smoothness parameter $\nu \geq 0.5$; $\lambda_C$ is the correlation length, $\sigma_{var}^2$ is the variance, $\Gamma$ is the gamma function and $K_\nu$ is the modified Bessel function of the second kind. The limiting case, i.e. $\nu \to \infty$, corresponds to the Gaussian covariance function $C_\infty (x,y) \ = \ \sigma_{var}^2 
\exp ( -|x - y|^2/\lambda_C^2 )$.

To sample from $\sigma_S$ we use the Karhunen-Lo{\`e}ve (KL) expansion of $\log \sigma_S$\,, i.e.,
\begin{equation}
\label{eq:KLexp}
\log \sigma_S (x,\omega) \ = \ \sum_{i=1}^{\infty} \sqrt{\xi_i} \ \eta_i(x) \ Z_i(\omega) \ ,
\end{equation}
where $Z_i \sim \mathcal{N}(0,1)$ i.i.d. Here  $\xi_i$ and $\eta_i$ are the eigenvalues and the $L^2(0,1)$-orthogonal eigenfunctions of the covariance integral operator associated with kernel given by 
the covariance function in \eqref{eq:covfunc}. In practice, the KL expansion needs to be truncated after a finite number of terms (here denoted $d$). The accuracy of this truncation depends on the decay of the eigenvalues \cite{Lord:14}. For $\nu < \infty$, this decay is algebraic and depends on the smoothness parameter $\nu$. In the Gaussian covariance case the decay is exponential.
Note that for the Mat{\'e}rn covariance with $\nu = 0.5$, the eigenvalues and eigenfunctions can be computed analytically \cite{Lord:14}. For other cases of $\nu$, we numerically compute the eigensystem using 
the Nystr{\"o}m method - see, for example,  \cite{EiErUl:07}.

The goal of stochastic uncertainty quantification is to understand how the randomness in $\sigma_S$ and $\sigma$ propagates to functionals of the scalar or angular flux. Such quantities of interest may be point values, integrals or norms of $\phi$ or $\psi$. They are random variables and the focus is on estimating their mean, variance or distribution.

\subsection{Discretisation}
\label{sec:gjquadrule}

For each realisation $\omega \in \Omega$, the stochastic 1D NTE  \eqref{eq:transport}, \eqref{eq:phi}, \eqref{eq:bc} is an integro-differential equation in two variables, space and angle. For ease of presentation, we suppress the dependency on $\omega \in \Omega$ for the moment. 

We use a  $2N$-point quadrature rule 
$\int_{-1}^1 f(\mu) d\mu \approx \sum_{|k|=1}^N w_k f(\mu_k)$ with nodes $\mu_k \in [-1,1]\backslash\{0\}$ and positive weights $w_k$ to discretise in angle, assuming the 
(anti-) symmetry properties $\mu_{-k} = -\mu_k$ and $w_{-k} = w_k$.  (In later sections, we construct such a rule by using $N$-point Gauss-Legendre rules on each of $[-1,0)$ and $(0,1]$.) 
   
To discretise in space, we introduce a mesh  $0 = x_0 < x_1 < \ldots < x_M = 1$ which is assumed to  resolve 
any  discontinuities in the cross-sections $\sigma, \sigma_S$ and is also     quasiuniform -  
i.e. the subinterval lengths $h_j : = x_j - x_{j-1}$ satisfy  
$\gamma h \leq h_j \leq h : = \max_{j = 1, \ldots  M} h_j , $ for some constant $\gamma >0$. 
Employing  a simple  Crank-Nicolson method for the transport part of \eqref{eq:transport}, \eqref{eq:phi} and combining it with the angular 
quadrature rule above we obtain the classical  \textit{diamond-differencing} scheme:
\begin{align}
\mu_k  &\frac{ \Psi_{k,j} - \Psi_{k,j-1}}{h_j} \ + \ \sigma_{j-1/2} \frac{\Psi_{k,j} + \Psi_{k,j-1}}{2}  \nonumber \\
&\qquad= \ \sigma_{S,j-1/2} \Phi_{j-1/2} \ + \ F_{j-1/2}   \ , \ \ \ \ \ j = 1,...,M, \ \ |k| = 1,\ldots,N, \label{eq:fulldisc} 
\end{align}
where
\begin{equation}
\label{eq:phi_fulldisc}
\Phi_{j-1/2}\ = \ \frac{1}{2} \sum_{|k| = 1}^N w_k \frac{\Psi_{k,j} \ + \ \Psi_{k,j-1}}{2} \ , \ \ j = 1,...,M \ .
\end{equation}
Here $\sigma_{j-1/2}$ denotes the value of $\sigma$ at the mid-point of the interval $I_j = (x_{j-1}, x_j)$, with the analogous  meaning for   
$\sigma_{S,j-1/2}$ and $F_{j-1/2}$. The notation reflects the fact that (in the next section) we will associate the unknowns $\Psi_{k,j}$ in \eqref{eq:fulldisc}  
with the nodal values $\psi_{k,h} (x_j)$ of continuous piecewise-linear functions  $\psi_{k,h} \approx \psi(\cdot, \mu_k)$. 

Finally, \eqref{eq:fulldisc} and \eqref{eq:phi_fulldisc} have  to be supplemented with the boundary conditions $\Psi_{k,0} = 0$, for $k >0$ and  
$\Psi_{k,M} = 0$, for $k <0$. If the right-hand side of \eqref{eq:fulldisc} were known, then   \eqref{eq:fulldisc} could be solved 
simply by sweeping from left to right (when $k>0$) and from right to left (when $k<0$). The appearance of $\Phi_{j-1/2}$ on the right-hand side means that \eqref{eq:fulldisc} and \eqref{eq:phi_fulldisc} consitute a 
coupled  system 
with solution  $(\Psi, \Phi)  \in \mathbb{R}^{2NM}\times \mathbb{R}^M$. It is helpful to think of  
$\Psi$ as being composed  of $2N$ subvectors $\Psi_k$, each  with $M$ entries $\Psi_{k,j}$, consisting of 
approximations to $\psi(x_j, \mu_k)$ with   $x_j$ ranging  over all free nodes.  

The coupled system \eqref{eq:fulldisc} and \eqref{eq:phi_fulldisc}
 can be written in matrix form as
\begin{equation}
\label{eq:fullsystem}
\begin{pmatrix}
 T & -\Sigma_S \\
 -P & I
\end{pmatrix}
\begin{pmatrix}
\Psi \\
\Phi
\end{pmatrix}
\ = \ \begin{pmatrix}
F \\
\vec{0}
\end{pmatrix} \ .
\end{equation}
Here, the vector $\Phi \in \mathbb{R}^M$ contains the approximations of the scalar flux at the $M$ midpoints of the spatial mesh. The matrix $T$ is a block diagonal $2NM \times 2NM$ matrix, representing the left hand side of \eqref{eq:fulldisc}. The $2N$ diagonal blocks of $T$, one per angle, are themselves bi-diagonal. The $2NM \times M$ matrix $\Sigma_S$ simply consists of $2N$ identical diagonal blocks, one per angle, representing the multiplication of $\Phi$ by $\sigma_S$ at the midpoints of the mesh. The $M \times 2NM$ matrix $P$ represents the right hand side of \eqref{eq:phi_fulldisc}, i.e. averaging at the midpoints and quadrature. The matrix $I$ denotes the $M \times M$ identity matrix. The vector $F \in \mathbb{R}^{2NM}$ contains $2N$ copies of the source term evaluated at the $M$ midpoints of the spatial mesh.

\subsection{Direct and Iterative Solvers}
\label{sec:solve}

We now wish to find the (approximate) fluxes in the linear system \eqref{eq:fullsystem}. We note that the matrix $T$ is invertible and has a useful sparsity structure that allows its inverse to be calculated in $\mathcal{O}(MN)$ operations. However, the bordered system \eqref{eq:fullsystem} is not as easy to invert, due to the presence of $\Sigma_S$ and $P$.

To exploit the sparsity of $T$,  we do block elimination on  \eqref{eq:fullsystem}  
 obtaining the Schur complement system for the scalar flux, i.e.,
\begin{equation}
\label{eq:intform}
\left( I - P T^{-1} \Sigma_S \right) \Phi \ = \ P T^{-1} F \ ,
\end{equation}
which now requires the inversion of a smaller (dense) matrix. Note that \eqref{eq:intform} is a finite-dimensional version of the reduction of the integro-differential equation \eqref{eq:transport}, \eqref{eq:phi} to the integral form of the NTE, 
see \eqref{eq:IE1}.
In this case, the two dominant computations with $\mathcal{O}(M^2 N)$ and $\mathcal{O}(M^3)$ operations respectively, are the triple matrix product $P T^{-1} \Sigma_S$ in the construction of the Schur complement and the $LU$ factorisation of the $M \times M$ matrix $\left( I - P T^{-1} \Sigma_S \right)$. This leads to a total 
\begin{equation}
\label{eq:A} 
\mbox{theoretical cost of the direct solver} \ \sim \ \mathcal{O}(M^2(M \ + \ N)) \ .
\end{equation} 
We note that for stability reasons (see \S \ref{sec:theo}, also  \cite{PiSc:83} in  a simpler context),
the number of spatial and angular points should be related. A suitable choice is $M\sim N$, leading to a cost of the direct solver of $\mathcal{O}(M^3)$ in general.


The second approach for solving \eqref{eq:fullsystem} is an iterative solver commonly referred to as \textit{source iteration}, cf. \cite{Bl:16}. The form of   \eqref{eq:intform} naturally suggests the iteration
\begin{equation}
\label{eq:iterative}
\Phi^{(k)} \ = \ P T^{-1} \left(  \Sigma_S \Phi^{(k-1)} \ + \ F \right) \ ,
\end{equation}
where $\Phi^{(k)}$ is the approximation at the $k$th iteration, with $\Phi^{(0)} = P T^{-1} F$. This can be seen as a discrete version of an iterative method for the integral equation \eqref{eq:IE1}.


In practice, we truncate after $K$ iterations. The dominant computations in the source iteration are the $K$ multiplications with $PT^{-1} \Sigma_S$. Exploiting the sparsity of all the matrices involved, these multiplications cost  $\mathcal{O}(M N)$ operations, leading to an overall
\begin{equation} \label{eq:B}
\mbox{theoretical cost of source iteration} \ \sim \ \mathcal{O}\left( M \, N \, K \right) \ .
\end{equation} 
Our numerical experiments in Section \ref{sec5} show that for $N = 2M$ the hidden constants in the two estimates \eqref{eq:A} and \eqref{eq:B}  are approximately the same. Hence, whether the iterative solver is faster than the direct solver depends on whether the number of iterations $K$ to obtain an accurate enough solution is smaller or larger than $M$.


There are sharp theoretical results on the convergence of source iteration for piecewise smooth cross-sections \cite[Thm 2.20]{Bl:16}. In particular, if $\phi^{(K)}(\omega)$ denotes the approximation to $\phi(\omega)$ after $K$ iterations, then 
\begin{equation}
\label{eq:sourceconv}
\bigg \| \sigma^{1/2} \left( \phi - \phi^{(K)} \right) \bigg \|_2 
\ \leq \ C'  \left(\eta \bigg\| \frac{\sigma_S}{\sigma} \bigg \|_{\infty}\right)^K\ , 
\end{equation}
for some constant $C'$ and $\eta \leq 1$. 
That is, the error decays geometrically with rate no slower than  the spatial  
maximum of $\sigma_S / \sigma$. This value depends on $\omega$ and will change pathwise. Using this result as a guide together with \eqref{eq:sigA}, we assume that the convergence of the $L^2$-error with respect to $K$ can be bounded by
\begin{equation}
\label{eq:theor_conv}
\| \phi \ - \ \phi^{(K)} \|_2 \ \le \ C  \bigg\| \frac{\sigma_S}{\sigma} \bigg \|_{\infty}^K \ ,
\end{equation}
for some constant $C$ that we will estimate numerically in Section \ref{sec5}.

\section{Summary of Theoretical Results}
\label{sec:theo} 

The  rigorous  analysis  of  UQ 
for  PDEs with random coefficients requires  estimates  
for the error when     
discretisations in physical space (e.g. by finite differences)
and probability space (e.g. by sampling techniques) are combined.    
The 
 physical error estimates typically need to be  probabilistic in form  
(e.g. estimates of expectation of the physical error). 
Such estimates  are quite well-developed  for   elliptic PDEs  - see for 
example \cite{ChSc:13} but this question is almost untouched for the transport equation (or more specifically the NTE). We outline 
here some   results which 
are  proved  in the forthcoming paper \cite{GrPaSc:17}. 
This paper proceeds by first  giving  an error analysis for \eqref{eq:transport_det}, \eqref{eq:phi_det} with variable cross-sections, 
which is explicit in $\sigma, \sigma_S$,  and then uses   this to derive probabilistic error 
estimates for the spatial discretisation \eqref{eq:fulldisc}, \eqref{eq:phi_fulldisc}. 
      
The numerical analysis of the NTE (and related integro-differential equation problems such as radiative transfer) dates back at least as far as the work of H.B. Keller \cite{Ke:60}. After a huge growth in the  
mathematics   literature in the 1970's and 1980's,  progress has been slower  
since. 
This is perhaps surprising,  since  
discontinuous Galerkin (DG) methods have enjoyed a massive recent renaissance 
and the solution of the neutron transport problem was one 
of the key  motivations behind the original introduction of DG \cite{ReHi:73}.  
Even today,  an error analysis of the NTE with variable 
(even deterministic) cross-sections 
(with  explicit dependence on  the data) is still not available, 
even for the model case of mono-energetic 1D slab geometry considered here. 

The fundamental paper on the  analysis of the discrete ordinates method  for the NTE is  
\cite{PiSc:83}. Here a  full analysis of the combined effect of angular and spatial discretisation is given  under the assumption that the cross-sections $\sigma$ and $\sigma_S$ in \eqref{eq:transport} are  constant. 
The delicate relation between spatial and angular discretisation parameters 
required to  achieve stability and convergence is described there. 
Later research e.g. \cite{As:98}, \cite{As:09}  
produced analogous results for models of increasing complexity  and 
in higher dimensions, but the proofs  were mostly confined to the case of cross-sections that are constant in space.
A   separate and related sequence of papers  (e.g. \cite{LaNe:82}, \cite{Vi:84},  and  \cite{AlViGa:89})
allow for variation in cross-sections, but error estimates explicit in this data  are not available there.

The results outlined  here are orientated to the case when $\sigma, \sigma_S$ 
have relatively rough fluctuations.    
As a precursor to attacking the random case, we first consider rough deterministic coefficients defined as follows. We assume that there is some partition of $[0,1]$ and that $\sigma, \sigma_S$ are $C^\eta$ functions on each subinterval of the partition (with $\eta \in   (0,1]$),  but that $\sigma, \sigma_S$ may be discontinuous across the break points.
We assume that the mesh $x_j$ introduced in \S  \ref{sec:gjquadrule} resolves these break points. (Here $C^\eta$ is the usual  H\"{o}lder space of 
index $\eta$ with norm $\Vert \cdot \Vert_\eta$.) We also assume that 
the source function  $f \in C^\eta$. 
    
When discussing the error when  \eqref{eq:fulldisc}, \eqref{eq:phi_fulldisc} is applied to \eqref{eq:transport_det}, 
\eqref{eq:phi_det}, it is useful to consider the ``pure transport'' problem: 
\begin{equation} \label{eq:tr}  \mu \frac{\rd u}{\rd x} + \sigma u = g, \quad  \text{with} \ u(0) = 0, \ \text{when} \ \mu > 0 \quad  \text{and}\  u(1) = 0 \ 
\text{when}\  \mu < 0,  \end{equation} 
 and with $g \in C$ a generic right-hand side (where $\mu$ is now a parameter). 
Application of the Crank-Nicolson method (as in \eqref{eq:fulldisc}) yields 
\begin{equation} \label{eq:trh}  \mu \left(\frac{U_j - U_{j-1}}{h_j}\right)   +  \sigma_{j-1/2}  \left(\frac{U_j + U_{j-1} }{2}\right)  = g_{j-1/2} \ , \ \mbox{ for } j = 1,...,M \ ,
\end{equation}
with analogous boundary conditions, where, for any continuous function $c$, we use  $c_{j-1/2}$ to denote  $c(x_{j-1/2})$. Letting $V^h$ denote the space of continuous piecewise linear functions with respect to the mesh $\{x_j\}$, \eqref{eq:trh} is equivalent to seeking a 
$u^h \in V^h$ (with nodal values $U_j$) such that 
$$\int_{I_j} \left(\mu \frac{\rd u^h }{\rd x } + {\widetilde{\sigma}} u^h \right)\ = \ \int_{I_j} \widetilde{g}\  , \quad j = 1, \ldots, M , \quad \text{where} \quad 
I_j = (x_{j-1}, x_j),  \  $$
and  $\widetilde{c}$ denotes the piecewise constant function with respect to the grid $\{x_j\}$ which interpolates 
$c$ at the mid-points of subintervals.

It is easy to show that both \eqref{eq:tr}  and \eqref{eq:trh} have unique solutions and we denote the respective  
solution operators by $\mathcal{S}_\mu$ and $\mathcal{S}^h_\mu$, i.e.  
$$ u = {\mathcal{S}}_\mu g \quad \text{and} \quad u^h = \mathcal{S}^h_\mu g \ .$$
Bearing in mind the angular averaging process in \eqref{eq:phi_det} and \eqref{eq:phi_fulldisc}, 
it is useful to then introduce the corresponding continuous and discrete  spatial operators:  
$$ (\mathcal{K} g)(x) := \frac12 \int_{-1}^{1} \left({\mathcal{S}}_\mu g\right)(x) \rd \mu, \quad  
\text{and} \quad  (\mathcal{K}^{h,N} g)(x)  \ = \ \frac12  \sum_{\vert k \vert = 1 }^N w_k (\mathcal{S}_{\mu_k}^h g)(x)\ .   $$
It is easy to see (and  well known classically - e.g. \cite{KaKe:77}) that   
$$ (\mathcal{K} g)(x) = \frac12 \int_{0}^1 E_1(\vert \tau(x,y) \vert) g(y) \rd y,   $$ 
\text{where} $E_1$ is the exponential integral and  the function $\tau(x,y) = \int_x^y \sigma$ is known as  the optical path. 
In fact (even when $\sigma$ is merely continuous), $\mathcal{K}$ is a compact Fredholm integral operator on a range of function spaces 
and $\mathcal{K}^{h,N}$ is a finite rank approximation to it. 
The study of these integral operators in the deterministic case is a classical topic, e.g. \cite{Sl:75}.
In the case of random $\sigma$, $\mathcal{K}$ is an integral operator with a random kernel which merits further investigation. 
   Returning  to \eqref{eq:transport_det},  \eqref{eq:phi_det},  we see readily that 
\begin{align}\label{eq:IE1} \psi(x, \mu) = \mathcal{S}_{\mu} (\sigma_S \phi + f)\ ,  \quad \text{so that} \quad \phi \ = \ 
\mathcal{K} (\sigma_S \phi + f ) . \end{align} 
Moreover  \eqref{eq:fulldisc} and \eqref{eq:phi_fulldisc} correspond to a discrete analogue of \eqref{eq:IE1} as follows. 
Introduce the  family of functions $\psi^{h,N}_k \in V^h$, $\vert k \vert = 1, \ldots, N$, by  requiring    $\psi^{h,N}_k$ 
to have nodal values $\Psi_{k,j}$.  
Then set\vspace{-2ex}  $$\phi^{h,N} :=  \frac12 \sum_{\vert k \vert =1}^N w_k \psi_k^{h,N}\ \in V^h  , $$
and it follows that \eqref{eq:fulldisc} and \eqref{eq:phi_fulldisc} may be rewritten (for each $j = 1,...,M$)
$$\int_{I_j} \left(\mu_k \frac{\rd \psi_k^{h,N} }{\rd x } + {\widetilde{\sigma}} \psi_k^{h,N}  \right)\ = \ 
\int_{I_j} \widetilde{g^{h,N}}\  , \quad  \text{where} \quad g^{h,N} = \sigma_S \phi^{h,N} + f\ . $$
and thus 
\begin{align}\label{eq:IE2} \psi_k^{h,N} = \mathcal{S}^{h}_{\mu_k} \left(\sigma_S \phi^{h,N} + f\right),   \quad \text{so that} \quad \phi^{h,N} = \mathcal{K}^{h,N} (\sigma_S \phi^{h,N} + f)
\ . \end{align} 

The numerical analysis of  \eqref{eq:fulldisc} and \eqref{eq:phi_fulldisc} is done by  
analysing (the second equation in) \eqref{eq:IE2} as an approximation of the second equation in  \eqref{eq:IE1}. This is studied in detail in \cite{PiSc:83} for constant $\sigma, \sigma_S$. 
In \cite{GrPaSc:17} we discuss the variable case,     
 obtaining  all estimates explicitly in  $\sigma, \sigma_S$. 
Elementary manipulation on \eqref{eq:IE1} and \eqref{eq:IE2} shows that
\begin{equation} 
\phi - \phi^{h,N} = ( I -\mathcal{K}^{h,N}\sigma_S )^{-1} (\mathcal{K} - \mathcal{K}^{h,N}) ( \sigma_S \phi + f) ,  \label{eq:IE3}
\end{equation} 
and so 
\begin{equation} 
\Vert \phi - \phi^{h,N}\Vert_\infty  \leq \Vert ( I -\mathcal{K}^{h,N}\sigma_S )^{-1}\Vert_\infty  \Vert (\mathcal{K} - \mathcal{K}^{h,N}) ( \sigma_S \phi + f)\Vert_\infty  .  \label{eq:IE4} 
\end{equation} 
The error analysis in \cite{GrPaSc:17} proceeds by estimating the two terms on the 
right-hand side of \eqref{eq:IE4}  separately. We summarise the results in the lemmas below. To avoid writing down the  technicalities   
(which will be given in detail in \cite{GrPaSc:17}), in the following results, we do not give the explicit dependence of the constants $C_i, \ i = 1,2,\ldots,$ on the cross sections $\sigma$ and $\sigma_S$. For simplicity we restrict our summary to the case when the right-hand side of \eqref{eq:trh} is the average of $g$ over $I_j$ (rather than the point value $g_{j-1/2}$). The actual scheme \eqref{eq:trh} is then analysed by a perturbation argument, see \cite{GrPaSc:17}.

\begin{lemma} Suppose    $N$    
is sufficiently large and  $h \log N$ is sufficiently small.  Then
\begin{equation} \label{eq:IE5a}   \Vert (I - \mathcal{K}^{h,N}\sigma_S)^{-1} \Vert_\infty \le \ C_1 \ ,  
 \end{equation} 
where $C_1$ depends on $\sigma$ and $\sigma_S$, but is independent of $h$ and $N$.
\end{lemma} 
\noindent 
{\em Sketch of proof}  \ \ The proof is obtained by first obtaining an estimate of the form \eqref{eq:IE5a} for the quantity 
$
\Vert ( I - \mathcal{K}\sigma_S)^{-1} \Vert_\infty $, and then showing that the perturbation 
  $\Vert \mathcal{K} - \mathcal{K}^{h,N} \Vert_\infty$ is small, when $N$ is sufficiently large and  $h \log N$ is sufficiently small.   (The constraint linking $h$ and $\log N$  arises because the transport equation  \eqref{eq:transport_det} 
has a singularity at $\mu = 0$.) The actual values of $h,N$ which are sufficient to  ensure  that the bound \eqref{eq:IE5a} holds  
depend on the cross-sections $\sigma$, $\sigma_S$.  

\begin{lemma} 
\begin{equation*}
\Vert (\mathcal{K}  - \mathcal{K}^{h,N})(\sigma_S \phi + f)   \Vert_\infty \ \le \  \left(C_2 \, h \log N +  
C_3\,  h^\eta   \ + \  C_4\,  \frac{1}{N} \right) \Vert f  \Vert_\eta \ , 
\end{equation*}
where $C_2, C_3, C_4$ depend again on $\sigma$ and $\sigma_S$, but are independent of $h, N$ and $f$.
\end{lemma}    
{\em Sketch of proof} \ \ 
Introducing the semidiscrete 
operator:
$$(\mathcal{K}^N g)(x) \ = \ \frac12 \sum_{\vert k \vert = 1}^N w_k (\mathcal{S}_{\mu_k} g)(x) $$
(corresponding to applying quadrature in angle but no  discretisation in space), we then write $ \mathcal{K} - \mathcal{K}^{h,N} = (\mathcal{K} - \mathcal{K}^{N}) + 
(\mathcal{K}^{N}  - \mathcal{K}^{h,N})$ and consider, separately, the semidiscrete error due to quadrature in angle: 
\begin{equation}\label{eq:IE6} 
(\mathcal{K} - \mathcal{K}^{N})(\sigma_S \phi + f)  = \frac{1}{2} \left( \int_{-1}^1 \psi(x,\mu) \rd \mu - \sum_{\vert k \vert = 1}^N w_k \psi(x, \mu_k) \right), 
\end{equation} 
and the spatial error  for a given $N$:  
\begin{equation} \label{eq:IE7} 
(\mathcal{K}^N - \mathcal{K}^{h,N} )(\sigma_S \phi + f) = \frac12 \sum_{\vert k \vert = 1} ^N w_k 
\left( \mathcal{S}_{\mu_k}  - \mathcal{S}^h_{\mu_k} \right) (\sigma_S \phi + f) . 
\end{equation}  

The estimate for \eqref{eq:IE6}   uses estimates for   the regularity of $\psi $ with respect to $\mu$ (which are explicit in the cross-sections), while \eqref{eq:IE7}   is estimated by proving  stability  of the Crank-Nicolson method and  a cross-section-explicit bound on $\Vert \phi\Vert_\eta$.


Putting together Lemmas 1 and 2, we obtain  the following.\\

\begin{theorem}\label{thm:error_det} 
Under the assumptions outlined above, 
\begin{equation*} 
\Vert \phi - \phi^{h,N} \Vert_\infty \ \le \  C_1 \left(C_2 \, h \log N +  
C_3 \,  h^\eta   \ + \  C_4 \,  \frac{1}{N} \right) \Vert f \Vert_\eta\;.
\end{equation*}
\end{theorem}     

Returning to the case when $\sigma, \sigma_S$ are random functions,  
this theorem provides  pathwise estimates for the error. 
In \cite{GrPaSc:17}, these are  turned into estimates in the corresponding Bochner space 
provided the coefficients 
$C_i$ are bounded in probability space. Whether this is the case depends on the choice of the random model for $\sigma, \sigma_S$. 

In particular, using the results in \cite[\S 2]{ChSc:13}, \cite{GrKu:15}, it can be shown that $C_i \in L^p(\Omega)$, for all $1 \le p < \infty$, for the specific choices of $\sigma$ and $\sigma_S$ in \S \ref{sec:ModelProblem}. Hence, we have: 



\begin{corollary} For all $1 \leq p < \infty$, 
$$\Vert \phi - \phi^{h,N} \Vert_{L^p(\Omega, L^\infty(0,1))} \ \le \  C
\left( h \log N +  
  h^\eta   \ + \    \frac{1}{N} \right) \Vert f \Vert_\eta \, ,$$
where $C$ is independent of $h, N$ and $f$.

\end{corollary}     

\section{Modern Variants of Monte Carlo}
\label{sec:MonteCarlo}

Let $Q(\omega) \in \mathbb{R}$ denote a functional of $\phi$ or $\psi$ representing a quantity of interest. 
We will focus on estimating $\bbE[Q]$, the  expected value of $Q$. Since we are not specific about what functionals we are considering, this includes also higher order moments or CDFs of quantities of interest. The expected value is a high-dimensional integral and the goal is to apply efficient quadrature methods in high dimensions. We consider Monte Carlo type sampling methods.

As outlined above, to obtain samples of $Q(\omega)$ the NTE has to be approximated numerically. First, the random scattering cross section $\sigma_S$ in \eqref{eq:transport} is sampled using the KL expansion of $\log \sigma_S$ in \eqref{eq:KLexp} truncated after $d$ terms. The stochastic dimension $d$ is chosen sufficiently high so that the truncation error is smaller than the other approximation errors. For each $n \in \mathbb{N}$, let $Z^{n} \in \mathbb{R}^d$ be a realisation of the multivariate Gaussian coefficient $Z := (Z_i)_{i=1,\ldots,d}$ in the KL expansion \eqref{eq:KLexp}. Also, denote by $Q_{h}(Z^{n})$ the approximation of the $n$th sample of $Q$ obtained numerically using a spatial grid with mesh size $h$ and $2N$ angular quadrature points. We assume throughout that $N \sim 1/h$, so there is a single discretisation parameter $h$. 

We will consider various unbiased, sample-based estimators $\widehat{Q}_h$ for the expected value $\bbE[Q]$ and we will quantify the accuracy of each estimator by its mean square error (MSE) $e(\widehat{Q}_h)^2$. Since $\widehat{Q}_h$ is assumed to be an unbiased estimate of $\bbE[ Q_h]$, i.e. $\bbE[ \widehat{Q}_h]=\bbE[ Q_h]$, the MSE can be expanded as
\begin{equation}
e(\widehat{Q}_h)^2 \ = \ \bbE\left[ (\widehat{Q}_h - \bbE[Q])^2 \right] \ = \ \left( \bbE\left[ Q - Q_h \right] \right)^2 \ + \ \mathbb{V}[\widehat{Q}_h] \ ,
\label{eq:MSE}
\end{equation}
i.e., the squared bias due to the numerical approximation plus the sampling (or quadrature) error $\mathbb{V}[\widehat{Q}_h] = \bbE[ (\widehat{Q}_h - \bbE[Q_h])^2 ]$.
In order to compare computational costs of the various methods we will consider their $\epsilon$-cost $\mathcal{C}_{\epsilon}$, that is, 
the number of floating point operations to achieve a MSE $e(\widehat{Q}_h)^2$ less than $\epsilon^2$. 

To bound the $\epsilon$-cost for each method, we make the following assumptions on the discretisation error and on the average cost to compute a sample from $Q_h$: 
\begin{align}
\Big|\bbE\left[ Q - Q_h \right]\Big| \ & = \ \mathcal{O}( h^\alpha) \ ,\label{eq:bias_assmc} \\
\bbE\left[\mathcal{C}(Q_h)\right] \ &\le \ \mathcal{O}(h^{-\gamma}) \ ,\label{eq:cost_assmc}
\end{align}
for some constants $\alpha, \gamma > 0$. We have seen in Section 2 that \eqref{eq:cost_assmc} holds with $\gamma$ between $2$ and $3$. The new theoretical results in Section 3 guarantee that \eqref{eq:bias_assmc} also holds for some $0 < \alpha \leq  1$. Whilst the results of Section 3 (and \cite{GrPaSc:17}) are shown to be sharp in some cases, the practically observed values for $\alpha$ in the numerical experiments here are significantly bigger, with values between 1.5 and 2. 

In recent years, many alternative methods for high-dimensional integrals have emerged that use tensor product deterministic quadrature rules combined with sparse grid techniques to reduce the computational cost \cite{XiKa:02, BaNoTe:07, NoTeWe:08, GuWeZh:14, AyEa:15, Fi:11, Gi6:13}. The efficiency of these approaches relies on high levels of smoothness of the parameter to output map and in general their cost may grow exponentially with the number of parameters (the \textit{curse of dimensionality}). Such methods are not competitive with Monte Carlo type methods for problems with low smoothness in the coefficients, where large numbers of parameters are needed to achieve a reasonable accuracy. For example, in our later numerical tests we will consider problems in up to 3600 stochastic dimensions.


However, standard Monte Carlo methods are notoriously slow to converge, requiring thousands or even millions of samples to achieve acceptable accuracies. In our application, where each sample involves the numerical solution of an integro-differential equation this very easily becomes intractable. The novel Monte Carlo approaches that we present here, aim to improve this situation in two complementary ways. Quasi-Monte Carlo methods reduce the number of samples to achieve a certain accuracy dramatically by using deterministic ideas to find well distributed samples in high dimensions. Multilevel methods use the available hierarchy of numerical approximations to our integro-differential equation to shift the bulk of the computations to cheap, inaccurate coarse models while providing the required accuracy with only a handful of expensive, accurate model solves.


\subsection{Standard Monte Carlo}
The (standard) Monte Carlo (MC) estimator for $\bbE[Q]$ is defined by
\begin{equation}
\widehat{Q}_h^{MC} \ := \ \frac{1}{N_{MC}} \sum_{n=1}^{N_{MC}} Q_h (Z^{n}) \ ,
\label{eq:estimatorMC}
\end{equation}
where $N_{MC}$ is the number of Monte Carlo points/samples $Z^{n} \sim \mathcal{N}(0,I)$. The sampling error of this estimator is $\mathbb{V}[\widehat{Q}_h^{MC}] = \mathbb{V}[Q_h]/N_{MC}$. 

A sufficient condition for the MSE to be less than $\epsilon^2$ is for both the squared bias and the sampling error in \eqref{eq:MSE} to be less than $\epsilon^2/2$. 
Due to assumption \eqref{eq:bias_assmc}, a sufficient condition for the squared bias to be less than $\epsilon^2 / 2$ is $h \sim \epsilon^{1/\alpha}$. Since $\mathbb{V}[Q_h]$ is bounded with respect to $h \to 0$, the sampling error of $\widehat{Q}_h^{MC}$ is less than $\epsilon^2 / 2$ for $N_{MC} \sim \epsilon^{-2}$. 
With these choices of $h$ and $N_{MC}$, it follows from Assumption \eqref{eq:cost_assmc} that the mean $\epsilon$-cost of the standard Monte Carlo estimator is 
\begin{equation}
\bbE\left[\mathcal{C}_{\epsilon}(\widehat{Q}_h^{MC})\right] \ = \ \bbE\left[ \sum_{n=1}^{N_{MC}}\mathcal{C}(Q_h(Z^{n}))\right] \ = \ N_{MC} \, \bbE\left[\mathcal{C}(Q_h)\right] \ = \ \mathcal{O}\left( \epsilon^{-2 - \frac{\gamma}{\alpha}}\right) \ .
\label{eq:epscost_MC}
\end{equation}
Our aim is to find alternative methods that have a lower $\epsilon$-cost.

\subsection{Quasi-Monte Carlo}
\label{sec:qmc}

The first approach to reduce the $\epsilon$-cost is based on using quasi-Monte Carlo (QMC) rules, which replace the random samples in \eqref{eq:estimatorMC} by carefully chosen deterministic samples and treat the expected value with respect to the $d$-dimensional Gaussian $Z$ in \eqref{eq:KLexp} as a high-dimensional integral with Gaussian measure. 

Initially interest in QMC points arose within number theory in the 1950's, and the theory is still at the heart of good QMC point construction today. Nowadays, the fast component-by-component construction (CBC) \cite{NuCo:06} provides a quick method for generating good QMC points, in very high-dimensions. Further information on the best choices of deterministic points and QMC theory  can be found in e.g. \cite{SlWo:98, DiPi:10, Ni:10, DiFr:13}. 

The choice of QMC points can be split into two categories; lattice rules and digital nets. We will only consider randomised rank-1 lattice rules here. 
In particular, given a suitable generating vector $z \in \mathbb{Z}^d$ and $R$ independent, uniformly distributed random shifts $(\Delta_r)_{r=1}^R$ in $[0,1]^d$, we construct $N_{QMC} =  R \, P$ lattice points in the unit cube $[0,1]^d$ using the simple formula
\[
v^{(n)} = \text{frac}\left( \frac{n z}{P} + \Delta_r \right), \qquad n=1,\ldots,P,\ \ r=1,\ldots,R
\]
where ``$\text{frac}$'' denotes the fractional part function applied componentwise and the number of random shifts $R$ is fixed and typically small e.g. $ R = 8, 16$.
To transform the lattice points $v^{n} \in [0,1]^d$ into ``samples'' $\widetilde{Z}^{n} \in \mathbb{R}^d$, $n=1,\ldots,N_{QMC}$, of the multivariate Gaussian coefficients $Z$ in the KL expansion \eqref{eq:KLexp} we apply the inverse cumulative normal distribution. See \cite{Gr:11} for details.

Finally, the QMC estimator is given by
$$
\widehat{Q}_h^{QMC} \ := \ \frac{1}{N_{QMC}} \sum_{n=1}^{N_{QMC}} Q_h(\widetilde{Z}^{n}) \ ,
$$
Note that this is essentially identical in its form to the standard MC estimator \eqref{eq:estimatorMC}, but crucially with deterministically chosen and then randomly shifted $\widetilde{Z}^{n}$.
The random shifts ensure that the estimator is unbiased, i.e. $\bbE[\widehat{Q}_h^{QMC}] = \bbE[Q_h]$.

The bias for this estimator is identical to the MC case, leading again to a choice of $h \sim \varepsilon^{1/\alpha}$ to obtain a MSE of $\varepsilon^2$. Here the MSE corresponds to the
mean square error of a randomised rank-1 lattice rule with $P$ points 
averaged over the shift $\Delta \sim \mathcal{U}([0,1]^d)$. In many cases, it can be shown that the quadrature error, i.e., the second term in \eqref{eq:MSE}, converges with $\mathcal{O}(N_{QMC}^{-1/2\lambda})$, with $\lambda \in (\frac{1}{2} ,1]$. That is, we can potentially achieve $\mathcal{O}(N_{QMC}^{-1})$ convergence for $\widehat{Q}_h^{QMC}$ as opposed to the $\mathcal{O}(N_{MC}^{-1/2})$ convergence for $\widehat{Q}_h^{MC}$. A rigorous proof of the rate of convergence requires detailed analysis of the quantity of interest (the integrand), in an appropriate weighted Sobolev space, e.g. \cite{GrKu:15}. Such an analysis is still an open question for this class of problems, and we do not attempt it here. Moreover, the generating vector $z$ does in theory have to be chosen problem specific. However, standard generating vectors, such as those available  at \cite{Kuo:QMClattice}, seem to also work well (and better than MC samples). Furthermore, we note the recent developments in ``higher-order nets'' \cite{GoDi:15,DiKuGi:14a}, which potentially increase the convergence of QMC methods to $\mathcal{O}(N_{QMC}^{-q})$, for $q \ge 2$.


Given the improved rate of convergence of the quadrature error and fixing the number of random shifts to $R = 8$, it suffices to choose $P \ \sim \ \epsilon^{-2 \lambda}$ for the quadrature error to be $\mathcal{O}(\varepsilon^2)$. Therefore it follows again from Assumption \eqref{eq:cost_assmc} that the $\epsilon$-cost of the QMC method satisfies
\begin{equation}
\bbE_\Delta\left[\mathcal{C}_{\epsilon}(\widehat{Q}^{QMC}) \right] \ = \ \mathcal{O}\left( \epsilon^{- 2\lambda -  \frac{\gamma}{\alpha}}\right) \ .
\label{eq:epscost_QMC}
\end{equation}
When $\lambda \to \frac{1}{2}$, this is essentially a reduction in the $\epsilon$-cost by a whole order of $\epsilon$. In the case of non-smooth random fields, we typically have $ \lambda \approx 1$ and the $\epsilon$-cost grows with the same rate as that of the standard MC method. However, in our experiments and in experiments for diffusion problems \cite{Gr:11}, the absolute cost is always reduced.

\subsection{Multilevel Methods}

The main issue with the above methods is the high cost for computing the samples $\{ Q_h(Z^{(n)})\}$, each requiring us to solve the NTE. The idea of the multilevel Monte Carlo (MLMC) method is to use a hierarchy of discrete models of increasing cost and  accuracy, corresponding to a sequence of decreasing discretisation parameters $h_0 > h_1 > ... > h_L = h$. Here, only the most accurate model on level $L$ is designed to give a bias of $\mathcal{O}(\epsilon)$ by choosing $h_L = h \sim \epsilon^{1/\alpha}$ as above.  The bias of the other models can be significantly higher.

MLMC methods were first proposed in an abstract way for high-dimensional quadrature by Heinrich \cite{He:01} and then popularised in the context of stochastic differential equations in mathematical finance by Giles \cite{Gi:08}. MLMC methods were first applied in uncertainty quantification in \cite{BaScZo:11,Cl:11}. The MLMC method has quickly gained popularity and has been further developed and applied in a variety of other problems. See \cite{Gi:15} for a comprehensive review. In particular, the multilevel approach is not restricted to standard MC estimators and can also be used in conjunction with QMC estimators \cite{GiWa:09,KuScSl:15,Kuo:15} or with stochastic collocation \cite{TeJa:15}. Here, we consider multilevel variants of standard MC and QMC.

MLMC methods exploit the linearity of the expectation, writing
$$
\bbE[Q_h] \ = \ \sum_{\ell=0}^L \bbE[Y_\ell]\ , \qquad \text{where} \ \ Y_\ell := Q_{h_\ell} - Q_{h_{\ell-1}} \ \ \text{and} \ \ Q_{h_{-1}} := 0.
$$
Each of the expected values on the right hand side is then estimated separately. In particular, in the case of a standard MC estimator with $N_\ell$ samples for the $\ell$th term, we obtain the MLMC estimator
\begin{equation}
\widehat{Q}_h^{MLMC} \ := \ \sum_{\ell = 0}^L \widehat{Y}^{MC}_\ell \ = \
\sum_{\ell=0}^L \frac{1}{N_\ell} \sum_{n=1}^{N_\ell} Y_\ell (Z^{\ell,n}) \ .
\label{eq:estimatorMLMC}
\end{equation}
Here, $\{Z^{\ell,n}\}_{n=1}^{N_\ell}$ denotes the set of i.i.d. samples on level $\ell$, chosen independently from the samples on the other levels.

The key idea in MLMC is to avoid estimating $\bbE[Q_h]$ directly. Instead, the expectation $\bbE[Y_0] = \bbE[Q_{h_0}]$ of a possibly strongly biased, but cheap approximation of $Q_h$ is estimated. The bias of this coarse model is then estimated by a sum of correction terms $\bbE[Y_\ell]$ using increasingly accurate and expensive models. Since the $Y_\ell$ represent small corrections between the coarse and fine models, it is reasonable to conjecture that  there exists $\beta > 0$ such that 
\begin{equation}
\label{eq:var_assmc}
\mathbb{V}[Y_\ell] \ = \ \mathcal{O}(h_\ell^{\beta})\ ,
\end{equation}
i.e., the variance of $Y_\ell$ decreases as $h_\ell \to 0$. This is verified for diffusion problems in \cite{ChSc:13}. Therefore the number of samples $N_\ell$ to achieve a prescribed accuracy on level $\ell$ can be gradually reduced, leading to a lower overall cost of the MLMC estimator. More specifically, we have the following cost savings:
\begin{itemize}
\item On the coarsest level, using \eqref{eq:cost_assmc}, the cost per sample is reduced from $\mathcal{O}(h^{-\gamma})$ to $\mathcal{O}(h_0^{-\gamma})$. Provided $\mathbb{V}[Q_{h_0}] \approx  \mathbb{V}[Q_{h}]$ and $h_0$ can be chosen independently of $\epsilon$, the cost of estimating $\bbE[Q_{h_0}]$ to an accuracy of $\varepsilon$ in \eqref{eq:estimatorMLMC} is reduced to $\mathcal{O}(\epsilon^{-2})$.
\item On the finer levels, the number of samples $N_\ell$ to estimate $\bbE[Y_\ell]$ to an accuracy of $\varepsilon$ in \eqref{eq:estimatorMLMC} is proportional to $\mathbb{V}[Y_\ell] \epsilon^{-2}$. Now, provided $\mathbb{V}[Y_\ell] = \mathcal{O}(h_\ell^\beta)$, for some $\beta > 0$, which is guaranteed if $Q_{h_\ell}$ converges almost surely to $Q$ pathwise, then we can reduce the number of samples as $h_\ell \to 0$. Depending on the actual values of $\alpha, \; \beta$ and $\gamma$, the cost to estimate $\bbE[Y_L]$ on the finest level can, in the best case, be reduced to $\mathcal{O}(\epsilon^{-\gamma/\alpha})$. 
\end{itemize}

The art of MLMC is to balance the number of samples across the levels to minimise the overall cost. This is a simple constrained optimisation problem to achieve $\mathbb{V} [\widehat{Q}_h^{MLMC}] \le \epsilon^2/2$. As shown in \cite{Gi:08}, using the technique of Lagrange Multipliers, the optimal number of samples on level $\ell$ is given by
\begin{equation}
N_\ell \ = \ \left\lceil 2 \epsilon^{-2} \left( \sum_{\ell=0}^L \sqrt{\mathbb{V}[Y_\ell]/\mathcal{C}_\ell} \right) \sqrt{\mathbb{V}[Y_\ell]\mathcal{C}_\ell} \right\rceil \ ,
\label{eq:samples_giles}
\end{equation} 
where 
$\mathcal{C}_\ell := \bbE\left[\mathcal{C}(Y_\ell)\right]$. In practice, it is necessary to estimate 
$\mathbb{V}[Y_\ell]$ and $\mathcal{C}_\ell$ in \eqref{eq:samples_giles} from the computed samples, updating $N_\ell$ as the simulation progresses.

Using these values of $N_\ell$ it is possible to establish the following theoretical complexity bound for MLMC \cite{Cl:11}.

\begin{theorem}
Let us assume that \eqref{eq:bias_assmc}, \eqref{eq:var_assmc} and \eqref{eq:cost_assmc} hold with $\alpha, \beta, \gamma > 0$.
Then, with $L \sim \log(\epsilon^{-1})$ and with the choice of $\{ N_\ell \}_{l=0}^L$ in \eqref{eq:samples_giles} we have
\begin{equation}
\bbE \left[\mathcal{C}_{\epsilon}(\widehat{Q}_{h_L}^{MLMC})\right] \ = \ \mathcal{O}\Big( \epsilon^{- 2 -  \max \left( 0 , \ \frac{\gamma - \beta }{\alpha} \right)}\Big)  \ .
\label{eq:epscost_MLMC}
\end{equation}
When $\beta = \gamma$, then there is an additional factor $\log(\epsilon^{-1})$.
\end{theorem}

Using lattice points $\widetilde{Z}^{\ell,n}$, as defined in Section \ref{sec:qmc}, instead of the random samples $Z^{\ell,n}$ we can in the same way define a multilevel quasi-Monte Carlo (MLQMC) estimator
\begin{equation}
\widehat{Q}_h^{MLQMC} \ := \ \sum_{\ell = 0}^L \widehat{Y}^{QMC}_\ell \ = \
\sum_{\ell=0}^L \frac{1}{\widetilde{N}_\ell} \sum_{n=1}^{\widetilde{N}_\ell} Y_\ell (\widetilde{Z}^{\ell,n}) \ .
\label{eq:estimatorMLQMC}
\end{equation}
The optimal values for $\widetilde{N}_\ell$ can be computed in a similar way to those in the MLMC method. However, they depend strongly on the rate of convergence of the lattice rule and in particular on the value of $\lambda$ which is difficult to estimate accurately. We will give a practically more useful approach below.

It is again possible to establish a theoretical complexity bound, cf. \cite{KuScSl:15, Kuo:15}.

\begin{theorem}
Let us assume that \eqref{eq:bias_assmc} and \eqref{eq:cost_assmc} hold with $\alpha, \gamma > 0$ and that there exists $\lambda \in (\frac{1}{2},1]$ and $\beta > 0$ such that 
\begin{equation}
\label{eq:var_assqmc}
\mathbb{V}_\Delta[\widehat{Y}^{QMC}_\ell] \ = \ \mathcal{O}\left(\widetilde{N}_\ell^{-1/\lambda} h_\ell^{\beta}\right)\ .
\end{equation} 
Let the number of random shifts on each level be fixed to $R$ and let $L \sim \log(\epsilon^{-1})$. Then, there exists a choice of $\{ N_\ell \}_{l=0}^L$ such that
\begin{equation}
\bbE_\Delta \left[\mathcal{C}_{\epsilon}(\widehat{Q}_{h_L}^{MLQMC})\right] \ = \ \mathcal{O}\Big(\epsilon^{-2 \lambda - \max \left( 0, \ \frac{\gamma - \beta \lambda}{\alpha} \right)}\Big) \ .
\label{eq:epscost_MLQMC}
\end{equation}
When $\beta\lambda = \gamma$, then there is an additional factor $\log(\epsilon^{-1})^{1+\lambda}$.
\end{theorem}
The convergence rate can be further improved by using higher order QMC rules \cite{DiKuGi:14}, but we will not consider this here.

It can be shown, for the theoretically optimal values of $N_\ell$, that there exists a constant $C$ such that
\begin{equation}
\label{eq:equilib}
\frac{\mathbb{V}_\Delta[\widehat{Y}^{QMC}_\ell]}{\mathcal{C}_\ell} \ = \ C\ ,
\end{equation}
independently of the level $\ell$ and of the value of $\lambda$ (cf. \cite[Sect.~3.3]{Kuo:15}). The same holds for MLMC. This leads to the following adaptive procedure to choose $N_\ell$ suggested in \cite{GiWa:09}, which we use in our numerical experiments below instead of \eqref{eq:samples_giles} . 

In particular, starting with an initial number of samples on all levels, we alternate the following two steps until $\mathbb{V} [\widehat{Q}_h^{MLMC}] \le \epsilon^2/2$: 
\begin{enumerate}
\item[(i)] Estimate $\mathcal{C}_\ell$ and  $\mathbb{V}_\Delta[\widehat{Y}^{QMC}_\ell]$ (resp.~$\mathbb{V}[\widehat{Y}^{MC}_\ell]$).
\item[(ii)] Compute 
$$
\ell^* = \argmax_{\ell=0}^{L}\ \left( \frac{\mathbb{V}_\Delta[\widehat{Y}^{QMC}_\ell]}{\mathcal{C}_\ell} \right)
$$
and double the number of samples on level $\ell^*$.  
\end{enumerate}
This procedure ensures that, on exit, \eqref{eq:equilib} is roughly satisfied and the numbers of samples across the levels $N_\ell$ are quasi-optimal. 

We use this adaptive procedure for both the MLMC and the MLQMC method. The lack of optimality typically has very little effect on the actual computational cost. Since the optimal formula \eqref{eq:samples_giles} for MLMC also depends on estimates of $\mathcal{C}_\ell$ and  $\mathbb{V}[Y_\ell]$, it sometimes even leads to a better performance. An additional benefit in the case of MLQMC is that the quadrature error in rank-1 lattice rules is typically lowest when the numbers of lattice points is a power of 2.

\section{Numerical Results}
\label{sec5}

We now present numerical results to confirm the gains that are possible with the novel multilevel and quasi-Monte Carlo method applied to our 1D NTE model \eqref{eq:transport_det}, \eqref{eq:phi_det}, \eqref{eq:bc}. We assume that the scattering cross-section $\sigma_S$ is a log-normal random field as described in Section \ref{subsec:UQ} and that the absorption cross section is constant, $\sigma_A \equiv \exp(0.25)$. We assume no fission, $\sigma_F \equiv 0$, and a constant source term $f = \exp(1)$.  We consider two cases, characterised by the choice of smoothness parameter $\nu$ in the Mat{\'e}rn covariance function \eqref{eq:covfunc}. For the first case, we choose $\nu = 0.5$. This corresponds to the exponential covariance and in the following is called the ``exponential field''. For the second case, denoted the ``Mat{\'e}rn field'', we choose $\nu = 1.5$. 
The correlation length and the variance are $\lambda_C = 1$ and $\sigma_{var}^2 = 1$, respectively.
The quantity of interest we consider is 
\begin{equation}
\label{eq:qoi}
Q(\omega) = \int_0^1 \phi(x,\omega) \rd x \ .
\end{equation}

For the discretisation, we choose a uniform spatial mesh with mesh width $h = 1/M$ and a quadrature rule (in angle) with $2N = 4M$ points.
The KL expansion of $\log(\sigma_S)$ in \eqref{eq:KLexp} is truncated after $d$ terms. We heuristically choose $d$ to ensure that the error due to this truncation is negligible compared to the discretisation error. In particular, we choose $d = 8h^{-1}$ for the Mat{\'e}rn field and $d = 225h^{-1/2}$ for the exponential field, leading to a maximum of 2048 and 3600 KL modes, respectively, for the finest spatial resolution in each case. Even for such large numbers of KL modes, the sampling cost does not dominate because the randomness only exists in the (one) spatial dimension.

We introduce a hierarchy of levels $\ell = 0,...,L$ corresponding to a sequence of discretisation parameters $h_\ell = 2^{-\ell} h_0$ with $h_0 = 1/4$, and approximate the quantity of interest in \eqref{eq:qoi} by
$$
Q_h(\omega) \ := \ \frac{1}{M} \sum_{j=1}^M \Phi_{j-1/2}(\omega) \ .
$$

To generate our QMC points we use an (extensible) randomised rank-1 lattice rule (as presented in Section \ref{sec:qmc}), with $R=8$ shifts. We use the generating vector \texttt{\small lattice-32001-1024-1048576.3600}, which is downloaded from \cite{Kuo:QMClattice}.

\subsection{A Hybrid Direct-Iterative Solver}

To compute samples of the neutron flux and thus of the quantity of interest, we propose a hybrid version of the direct and the iterative solver for the Schur complement system \eqref{eq:intform}  described in Section \ref{sec:solve}. 

The cost of the iterative solver depends on the number $K$ of iterations that we take. For each $\omega$, 
we aim to choose $K$ such that the $L_2$-error $\| \phi(\omega) - \phi^{(K)}(\omega) \|_2$ is less than $\epsilon$. To estimate $K$ 
we fix $h=1/1024$ and $d=3600$ and use the direct solver to compute  $\phi_h$ for each sample $\omega$. Let $\rho(\omega) := \| \sigma_S(\cdot,\omega)/\sigma(\cdot,\omega) \|_{\infty}$.  For a sufficiently large number of samples,  we then evaluate 
\[
\frac{\log \Big(\big\| \phi_h(\omega) \ - \ \phi_h^{(K)}(\omega)  \big\|_2\Big)}{K \, \log \big( \rho(\omega)\big)} 
\]
and find that this quotient is less than $\log(0.5)$ in more than 99\% of the cases, for $K=1,\ldots,150$, so that we can choose $C=0.5$ in \eqref{eq:theor_conv}. We repeat the experiment also for larger values of $h$ and smaller values of $d$ to verify that this bound holds in at least 99\% of the cases independently of the discretisation parameter $h$ and of the truncation dimensions $d$. 

Hence, a sufficient, a priori condition to achieve $\| \phi_h(\omega) - \phi_h^{(K)}(\omega) \|_2 < \epsilon$ in at least 99\% of the cases is 
\begin{equation}
\label{eq:chooseK}
K \ = \  K( \epsilon, \omega) \ = \ \max \bigg \{ \ 1 , \
\bigg \lceil \frac{\log\left( 2\epsilon \right)}{\log \big(\rho(\omega)\big)} \bigg \rceil \ \bigg \} \ ,
\end{equation}
where $\lceil \cdot \rceil$ denotes the ceiling function.
It is important to note that $K$ is no longer a deterministic parameter for the solver (like $M$ or $N$). Instead, $K$ is a random variable that depends on the particular realisation of $\sigma_S$. It follows from \eqref{eq:chooseK}, using the results in \cite[\S 2]{ChSc:13}, \cite{GrKu:15} as in Section \ref{sec:theo}, that  $\bbE[K(\epsilon,\cdot)] = \mathcal{O}( \log(\epsilon) )$ and $\mathbb{V}[K(\epsilon,\cdot)] = \mathcal{O}\left( \log(\epsilon)^2 \right)$, with more variability in the case of the exponential field. 

Recall from \eqref{eq:A} and \eqref{eq:B} that, in the case of $N=2M$, the costs for the direct and iterative solvers are $C_1 M^3$ and $C_2 KM^2$, respectively. In our numerical experiments, we found that in fact $C_1 \approx C_2$, for this particular relationship between  $M$ and $N$. 
This motivates a third ``hybrid'' solver, presented in Algorithm~\ref{alg:hybrid}, where the iterative solver is chosen when $K(\omega) < M$ and the direct solver when $K(\omega) \geq M$. This allows us to use the optimal solver for each particular sample.
 
\begin{figure}[t]
\begin{center}
\includegraphics[width=0.485\textwidth]{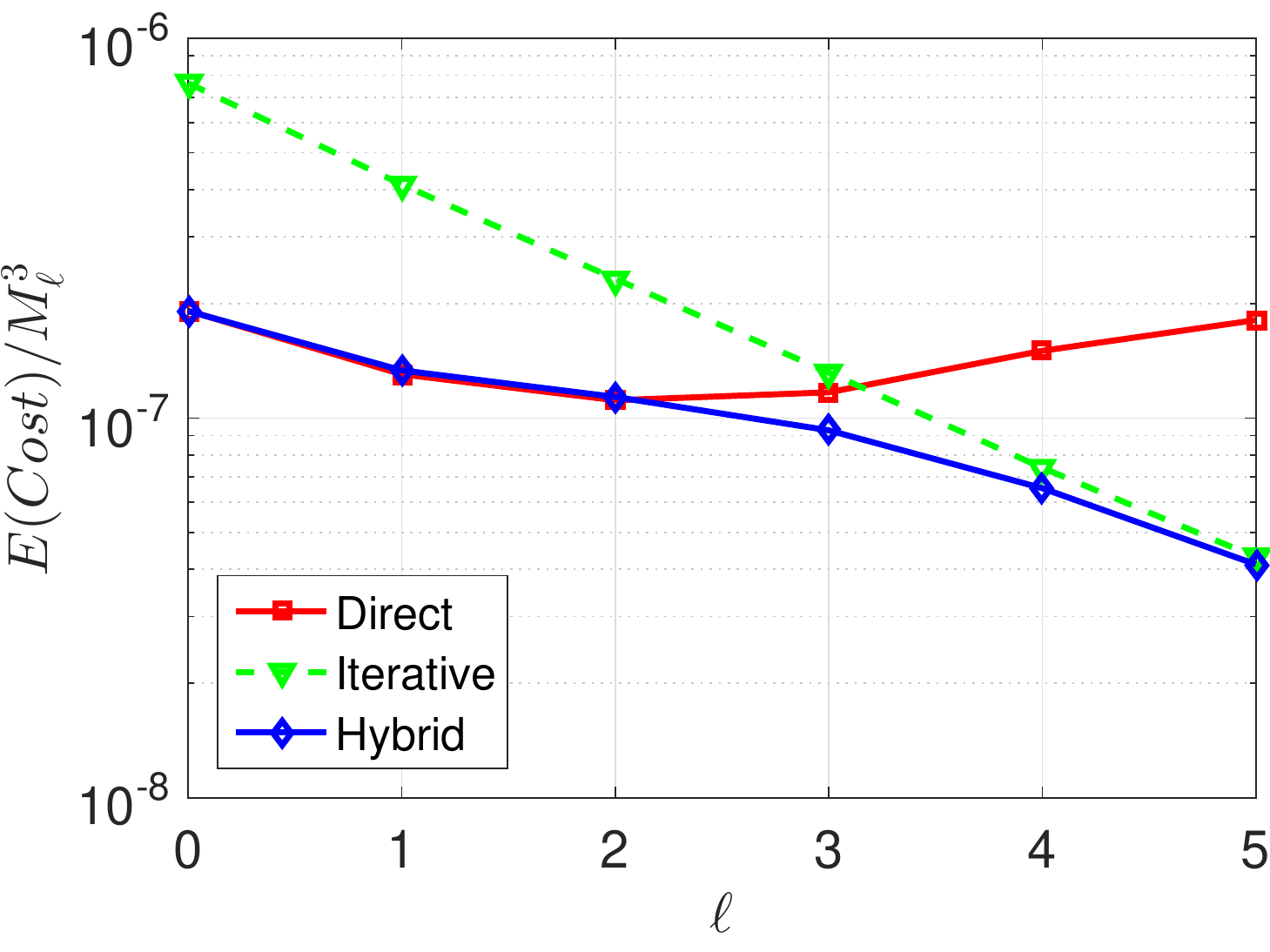}\hspace{2ex}\includegraphics[width=0.485\textwidth]{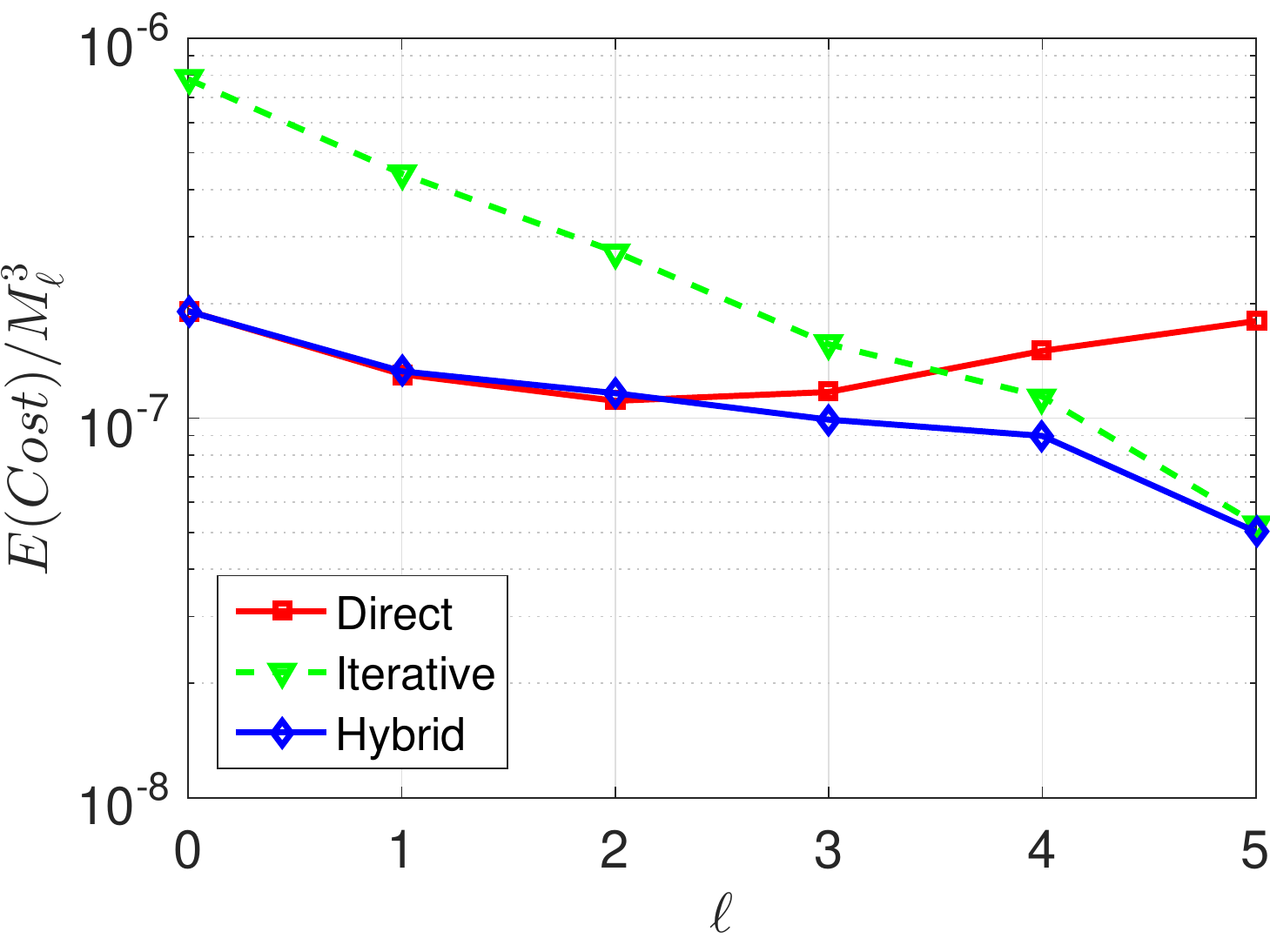}
\caption{Comparison of the average costs of the solvers (actual timings in seconds divided by $M_\ell^3$) for the Mat{\'e}rn field (left) and for the exponential field (right). \label{fig:scaledcost}}
\end{center}
\end{figure}

We finish this section with a study of timings in seconds (here referred to as the cost) of the three solvers. In Fig.~\ref{fig:scaledcost}, we plot the average cost (over $2^{14}$ samples) divided by $M_\ell^3$, against the level parameter $\ell$. We observe that, as expected, the (scaled) expected cost of the direct solver is almost constant and the iterative solver is more efficient for larger values of $M_\ell$. Over the range of values of $M_\ell$ considered in our experiments, a best fit for the rate of growth of the cost with respect to the discretisation parameter $h_\ell$ in \eqref{eq:cost_assmc} is $\gamma \approx 2.2$, for both fields. Thus our solver has a practical complexity of $\mathcal{O}(n^{1.1})$, where  $ n \sim M^2$ is the total number of degrees of freedom in the system.

\begin{algorithm}
\caption{ Hybrid direct-iterative solver of \eqref{eq:intform}, for one realisation}.          
\label{alg:hybrid} 
\begin{algorithmic}                    
    \REQUIRE Given $\sigma_S$, $\sigma$ and a desired accuracy $\epsilon$
    \STATE $K = \bigg \lceil \ log(2 \epsilon) \ / \ log(\rho) \ \bigg \rceil$
    \IF{$K < M$}
        \STATE Solve using $K$ source iterations
    \ELSE
        \STATE Solve using the direct method
    \ENDIF
\end{algorithmic}
\end{algorithm}

\vspace{0.2cm}

\subsection{A-Priori Error Estimates}
\label{sec:apriori}

Studying the complexity theorems of Section \ref{sec:MonteCarlo}, we can see that the effectiveness of the various Monte Carlo methods depends on the parameters $\alpha$, $\beta$, $\gamma$ and $\lambda$ in \eqref{eq:bias_assmc}, \eqref{eq:cost_assmc}, \eqref{eq:var_assmc} and \eqref{eq:var_assqmc}. In this section, we will (numerically) estimate these parameters in order to estimate the theoretical computational cost for each approach.

\begin{figure}[t]
\begin{center}
\includegraphics[width=0.485\textwidth]{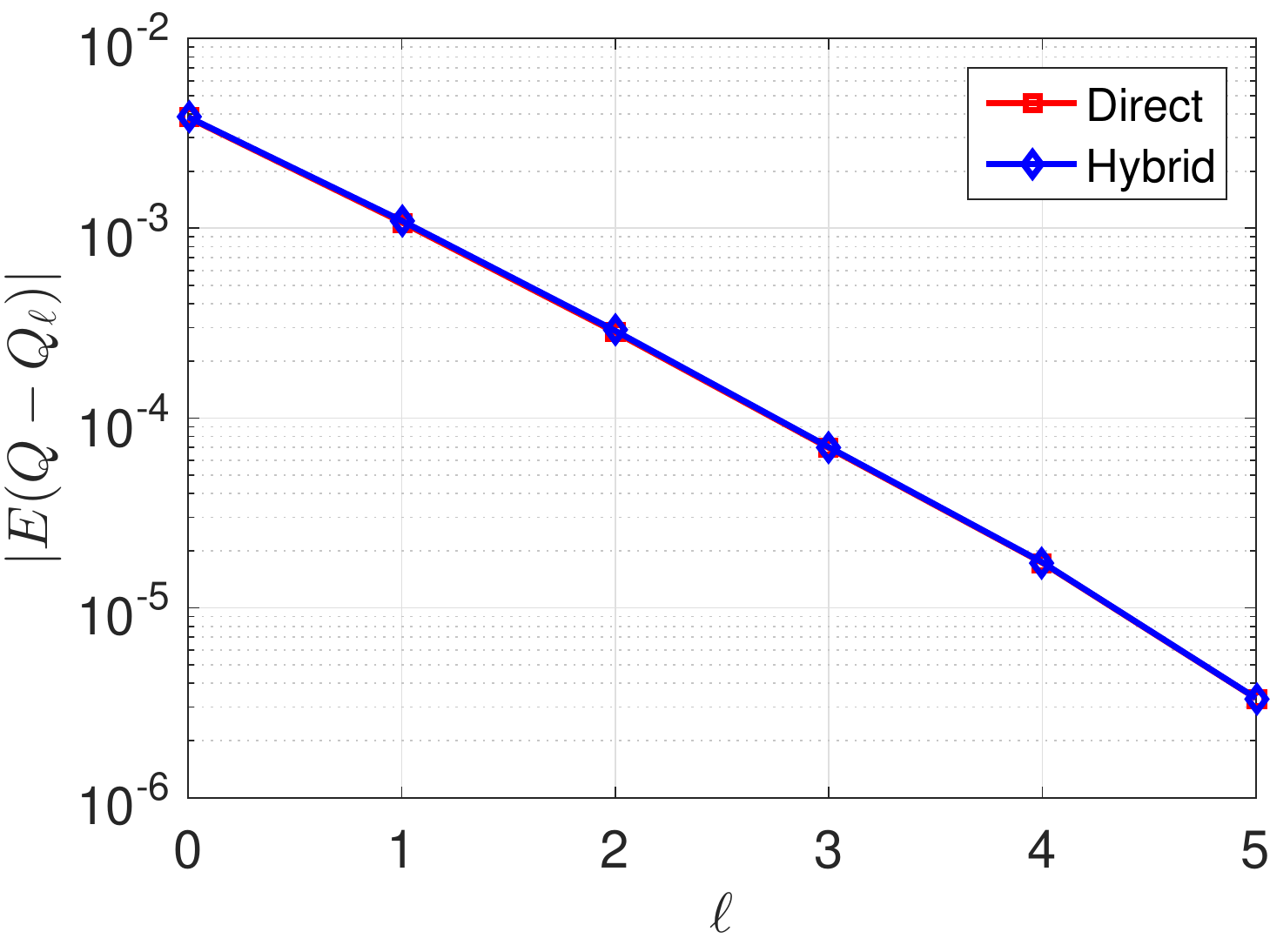}\hspace{2ex}\includegraphics[width=0.485\textwidth]{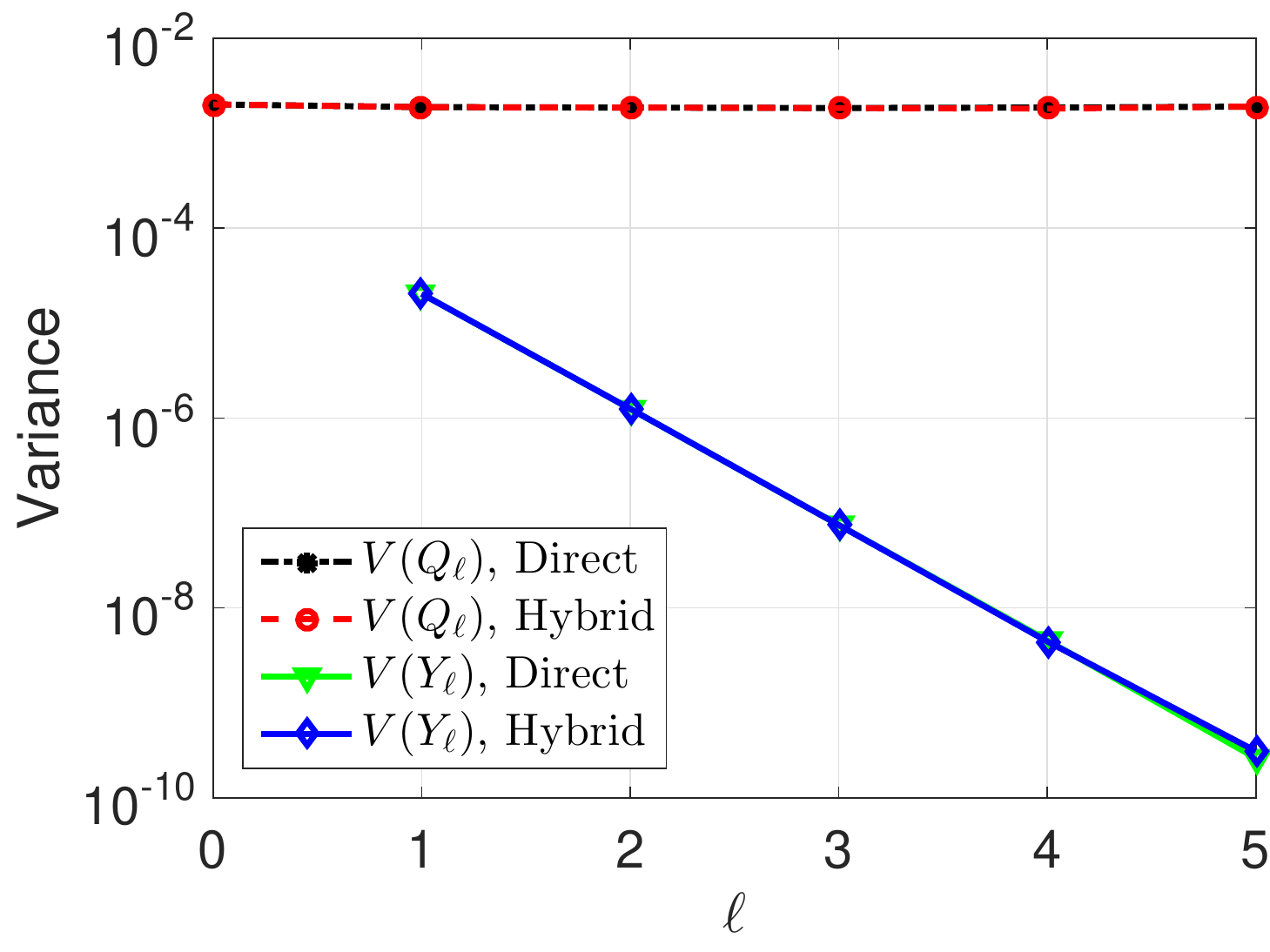}
\caption{Estimates of the bias due to discretisation errors (left) and of the variances of $Q_{h_\ell}$ and $Y_\ell$  (right), in the case of the Mat{\'e}rn field. \label{fig:biasvarcomp}}
\end{center}
\end{figure}

We have already seen that $\gamma \approx 2.2$ for the hybrid solver. In Fig.~\ref{fig:biasvarcomp}, we present estimates of the bias $\bbE[ Q - Q_{h_\ell}]$, as well as of the variances of $Q_{h_\ell}$ and of $Y_\ell$,  computed via sample means and sample variances over a sufficiently large set of samples.
We only explicitly show the curves for the Mat{\'e}rn field. The curves for the exponential field look similar. From these plots, we can estimate $\alpha \approx 1.9$ and $\beta \approx 4.1$, for the Mat{\'e}rn field, and $\alpha \approx 1.7$ and $\beta \approx 1.9$, for the exponential field. 

\begin{figure}[t]
\begin{center}
\includegraphics[width=0.485\textwidth]{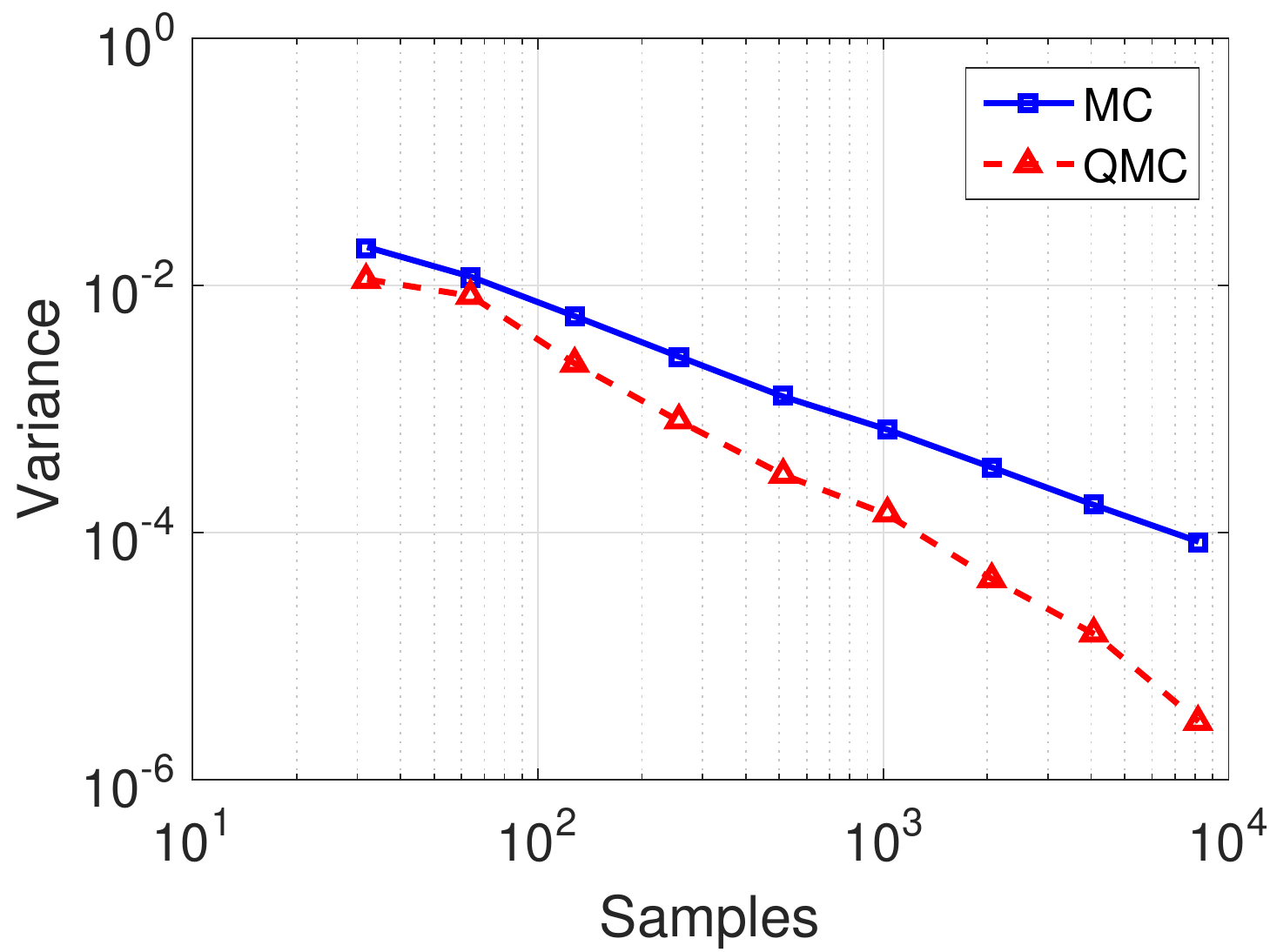}\hspace{2ex}\includegraphics[width=0.485\textwidth]{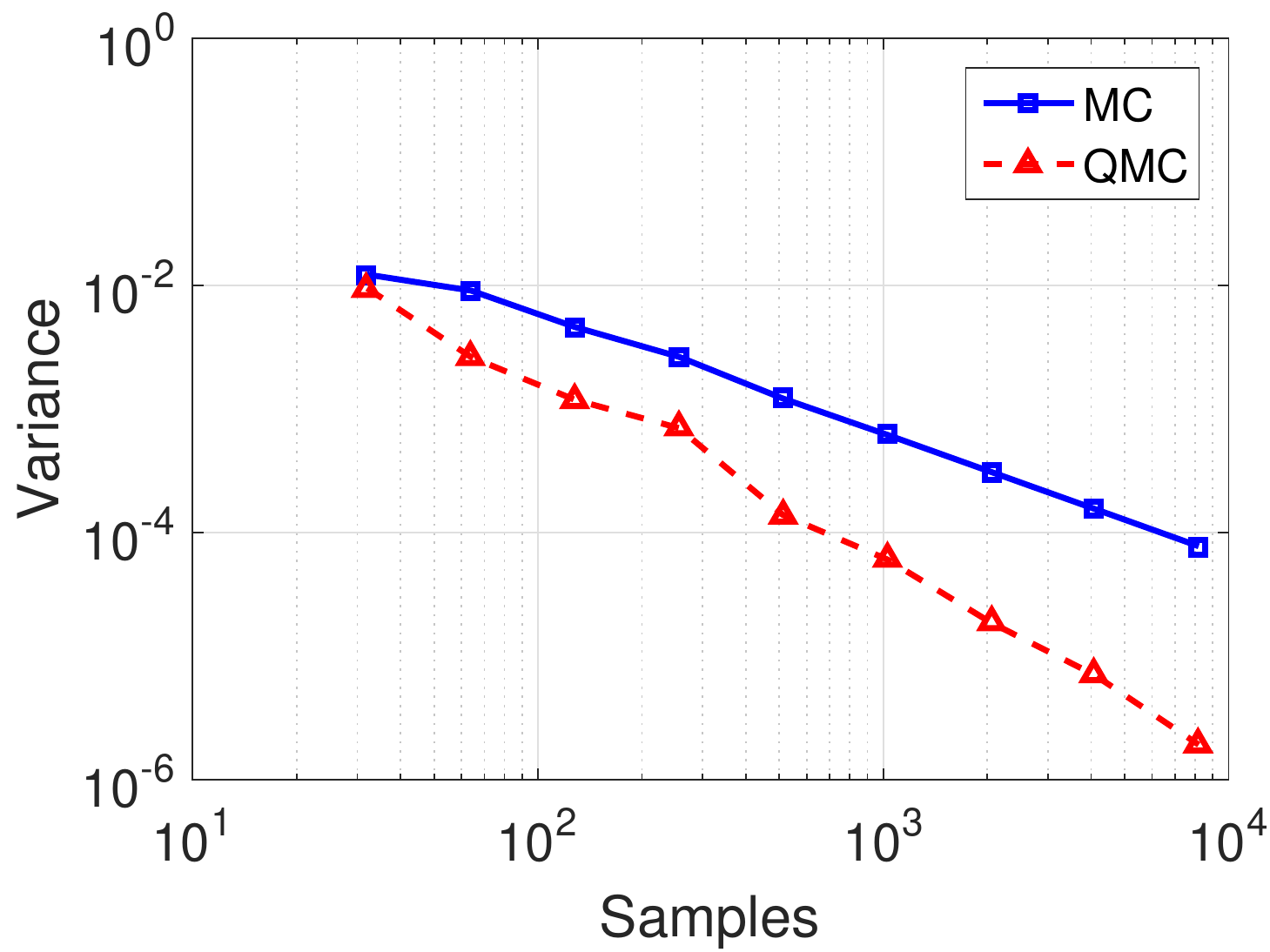}
\caption{Convergence of standard Monte Carlo and quasi-Monte Carlo estimators: Mat{\'e}rn field (left) and exponential field (right). \label{fig:mcvqmc}}
\end{center}
\end{figure}

To estimate $\lambda$ in \eqref{eq:var_assqmc}, we need to study the convergence rate of the QMC method with respect to the number of samples $N_{QMC}$. This study is illustrated  in Fig.~\ref{fig:mcvqmc}. As expected, the variance of the standard MC estimator converges with $\mathcal{O}(N_{MC}^{-1})$. On the other hand, we observe that the variance of the QMC estimator converges approximately with $\mathcal{O}(N_{QMC}^{-1.6})$ and $\mathcal{O}(N_{QMC}^{-1.4})$ (or $\lambda = 0.62$ and $\lambda = 0.71$) for the Mat{\'e}rn field and for the exponential field, respectively.

We summarise all the estimated rates in Table \ref{table:MCrates}.
\begin{table}[t]
\centering
\begin{tabular}{c|cccc}
 & $\alpha$ & $\beta$ & $\gamma$ & $\lambda$ \\\hline
Mat{\'e}rn field & \ 1.9 \ & \ 4.1 \ & \ 2.2 \ & \ 0.62 \\
Exponential field \ & \ 1.7 \ & \ 1.9 \ & \ 2.2 \ & \ 0.71 
 \end{tabular}
\caption{\label{table:MCrates} Summary of estimated rates  in \eqref{eq:bias_assmc}, \eqref{eq:cost_assmc}, \eqref{eq:var_assmc} and \eqref{eq:var_assqmc}.}
\end{table}

\subsection{Complexity Comparison of Monte Carlo Variants}

\begin{figure}[t]
\begin{center}
\includegraphics[width=0.485\textwidth]{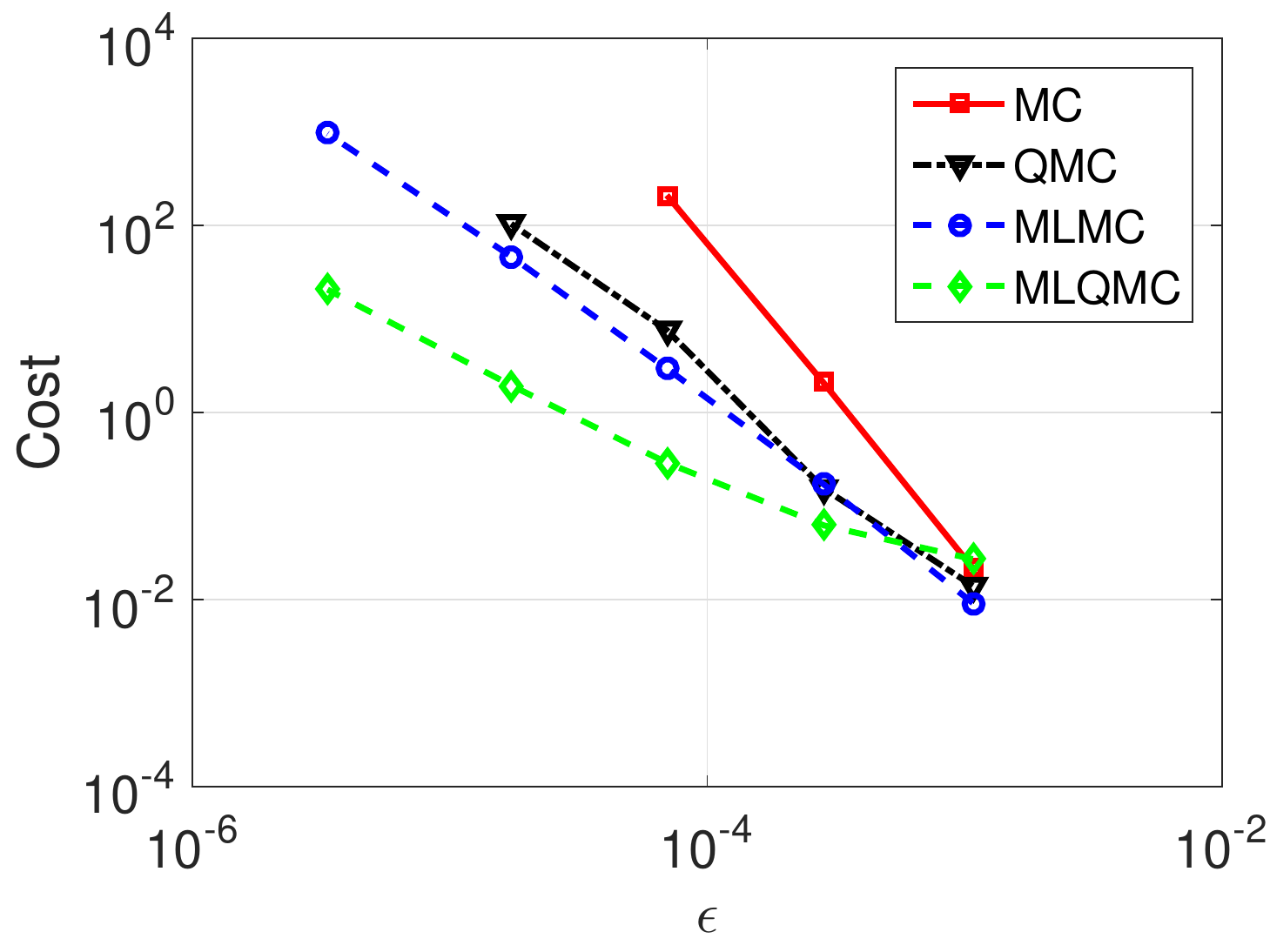}\hspace{2ex}\includegraphics[width=0.485\textwidth]{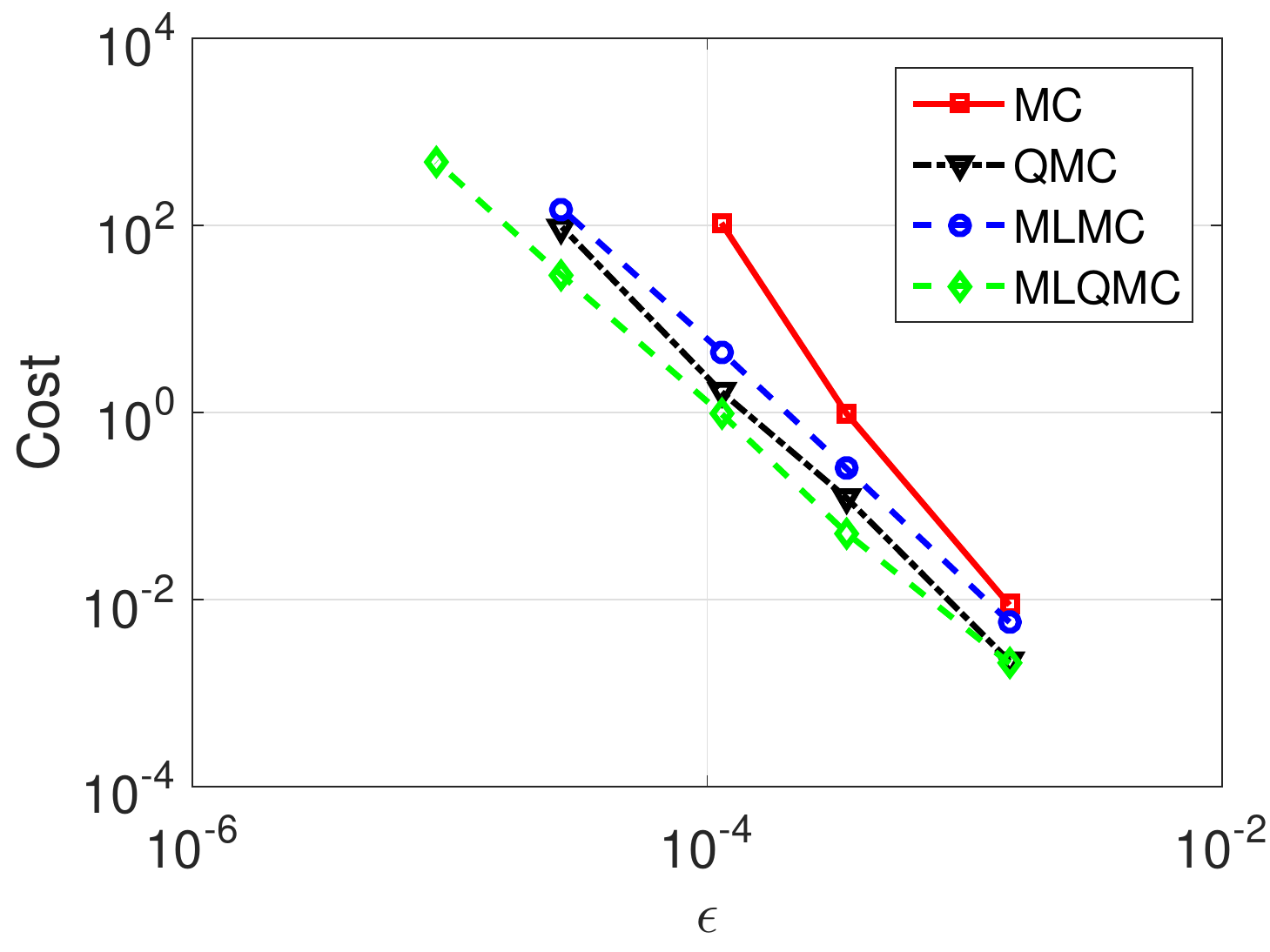}
\caption{Actual cost plotted against estimated bias on level $L$ for standard Monte Carlo, QMC, multilevel MC and multilevel QMC: Mat{\'e}rn field (left) and exponential field (right). \label{fig:epscost_hybrid}}
\end{center}
\end{figure}

For a fair comparison of the complexity of the various Monte Carlo estimators, we now use the a priori bias estimates in Section \ref{sec:apriori} to choose a suitable tolerance $\epsilon_L$ for each choice of $h = h_L$. Let $\tau_\ell$ be the estimated bias on level $\ell$. Then, for each $L=2,\ldots,6$, we choose $h = h_L$ and $\epsilon_L := \sqrt{2} \, \tau_L$, and we plot in Fig.~\ref{fig:epscost_hybrid} the actual cost of each of the estimators described in Section \ref{sec:MonteCarlo} against the estimated bias on level $L$. The numbers of samples for each of the estimators are chosen such that $\mathbb{V}[\widehat{Q}_h] \le \epsilon_\ell^2/2$. The coarsest mesh size in the multilevel methods is always $h_0 = 1/4$. We can clearly see the benefits of the QMC sampling rule and of the multilevel variance reduction, and the excellent performance of the multilevel QMC estimator confirms that the two improvements are indeed complementary. As expected, the gains are more pronounced for the smoother (Mat\'ern) field.

We finish by comparing the actual, observed $\epsilon$-cost of each of the methods with the $\epsilon$-cost predicted theoretically using the estimates for $\alpha$, $\beta$, $\gamma$ and $\lambda$ in Section~\ref{sec:apriori}. Assuming a growth of the $\epsilon$-cost proportional to $\epsilon^{-r}$, for some $r > 0$, we compare in Table \ref{table:MCepscost_hybrid} estimated and actual rates $r$ for all the estimators. Some of the estimated rates in Section~\ref{sec:apriori} are fairly crude, so the good agreement between estimated and actual rates is quite impressive.

\begin{table}[t]
\centering
\begin{tabular}{l||cc|cc|cc|cc}
             & \multicolumn{2}{c}{\textbf{MC}} & \multicolumn{2}{|c}{\textbf{QMC}} & \multicolumn{2}{|c}{\textbf{MLMC}} & \multicolumn{2}{|c}{\textbf{MLQMC}} \\
Field        &  Estimated          & \ Actual  \       &  Estimated          & \ Actual     \     & Estimated    &  \ Actual    \      & Estimated          & \ Actual  \       \\ \hline \hline
Mat{\'e}rn \ & 3.2            & 3.4            & 2.4             & 2.7            & 2.0             & 2.1             & 1.2              & 1.5             \\
Exponential \ & 3.3            & 3.6            & 2.7             & 2.4            & 2.2             & 2.5             & 1.9              & 1.9            
\end{tabular}
\caption{Comparison of the estimated theoretical and actual computational $\epsilon$-cost rates, for different Monte Carlo methods, using the hybrid solver.}
\label{table:MCepscost_hybrid}
\end{table}

\section{Conclusions}

To summarise, we have presented an overview of novel error estimates for the 1D slab geometry simplification of the Neutron Transport Equation, with spatially varying and random cross-sections. In particular, we consider the discrete ordinates method with Gauss quadrature for the discretisation in angle, and a diamond differencing scheme on a quasi-uniform grid in space. We represent the spatial uncertainties in the cross-sections by log-normal random fields with Mat{\'e}rn covariances, including cases of low smoothness. These  error estimates are the first of this kind. They allow us to satisfy key assumptions for the variance reduction in multilevel Monte Carlo methods. 

We then use a variety of recent developments in Monte Carlo methods to study the propagation of the uncertainty in the cross-sections, through to a linear functional of the scalar flux. We find that the Multilevel Quasi Monte Carlo method gives us significant gains over the standard Monte Carlo method. These gains can be as large as almost two orders of magnitude in the computational $\epsilon$-cost for $\epsilon = 10^{-4}$.

As part of the new developments, we present a hybrid solver, which automatically switches between a direct or iterative method, depending on the rate of convergence of the iterative solver which varies from sample to sample. Numerically, we observe that the hybrid solver is almost an order of magnitude cheaper than the direct solver on the finest mesh, on the other hand the direct solver is  
almost an order of magnitude cheaper than the iterative solver on the coarsest mesh we considered.

We conclude that modern variants of Monte Carlo based sampling methods are extremely useful for the problem of Uncertainty Quantification in Neutron Transport. This is particularly the case when the random fields are non-smooth and a large number of stochastic variables are required for accurate modelling.

\begin{acknowledgement}

We thank EPSRC and AMEC Foster Wheeler  for financial support
for this project and we particularly thank Professor Paul Smith (AMECFW) for many helpful
discussions. Matthew Parkinson is supported by the
EPSRC Centre for Doctoral Training in Statistical Applied Mathematics at
Bath (SAMBa), under project EP/L015684/1. This research made use of
the Balena High Performance Computing (HPC) Service
at the University of Bath.
\end{acknowledgement}

%
%


\bibliographystyle{plain}

\end{document}